\documentclass{amsart} 
\usepackage{amssymb, amscd}

\def\BB{{\mathbb B}}

\def\CC{{\mathbb C}}
\def\DD{{\mathbb D}}

\def\HH{{\mathbb H}}

\def\LL{{\mathbb L}}
 
\def\PP{{\mathbb P}}
\def\QQ{{\mathbb Q}} 
\def\RR{{\mathbb R}} 
\def\ZZ{{\mathbb Z}} 
\def\G{{\Gamma}}
\def\g{{\gamma}}

\def\an{{\rm an}}

\def\red{{\rm red}}

\def\st{{\rm st}}
\def\ss{{\rm ss}}

\def\tor{{tor}}
\def\bb{{bb}}

\def\bs{\backslash}
\def\bss{{\bs\!\bs}}
\def\dot{{\bullet}}
\def\pt{{\bullet}}

\def\Hcal{{\mathcal H}} 
\def\Ical{{\mathcal I}}

\def\Lcal{{\mathcal L}}

\def\Mcal{{\mathcal M}}

\def\Ocal{{\mathcal O}}

\def\Tcal{{\mathcal T}}

\def\Hbar{\overline{H}}

\def\la{\langle}
\def\ra{\rangle}

\newcommand\codim{\operatorname{codim}} 
 
\newcommand\Hom{\operatorname{Hom}}

\newcommand\lie{\operatorname{Lie}}
\newcommand\pic{\operatorname{Pic}}

\newcommand\GL{\operatorname{GL}}
\newcommand\PO{\operatorname{PO}}
\newcommand\SL{\operatorname{SL}}

\newcommand\U{\operatorname{U}}
\newcommand\SU{\operatorname{SU}}

\newtheorem{theorem}{Theorem}[section]
\newtheorem{lemma}[theorem]{Lemma}
\newtheorem{proposition}[theorem]{Proposition}
\newtheorem{corollary}[theorem]{Corollary}

\theoremstyle{definition}
\newtheorem{definition}[theorem]{Definition}
\newtheorem{convention}[theorem]{Convention}
\newtheorem{example}[theorem]{Example}

\theoremstyle{remark}

\newtheorem{question}[theorem]{Question}

\begin{document}

\title{Compactifications defined by arrangements I: the
ball quotient case}
\author{Eduard Looijenga}

\address{Faculteit Wiskunde en Informatica, 
Universiteit Utrecht, 
Postbus 80.010, NL-3508 TA Utrecht, Nederland}
\email{looijeng@math.uu.nl}



\date{\today}

\begin{abstract}
We define a natural compactification 
of an arrangement complement in a ball quotient. 
We show that when this complement has a moduli space
interpretation, then  this compactification is often one 
that appears naturally by means of geometric invariant theory.
We illustrate this with the moduli spaces of smooth quartic curves and
rational elliptic surfaces.
\end{abstract}
\maketitle

\section*{Introduction}
We wish to compare two compactifications 
of an algebro-geometric nature: those obtained by means of \emph{geometric
invariant theory (GIT)} and the ones defined 
by the \emph{Baily-Borel theory (BB)}, in situations
were both are naturally defined. The examples that we have in mind are 
moduli spaces of varieties for which a period mapping identifies 
that space with an open-dense
subvariety of an arithmetic quotient of a type IV domain or a complex ball.
This includes  moduli spaces of (perhaps multiply) polarized 
$K3$-surfaces, Enriques surfaces and Del 
Pezzo surfaces, but excludes for instance the moduli spaces of principally 
polarized abelian varieties of dimension at least three. In these examples, the 
GIT-boundary parametrizes degenerate forms of varieties, so that  
this boundary is  naturally stratified with strata parametrizing
degenerate varieties of the same type. The strata and their incidence
relations have been computed in a number of cases, sometimes after 
hard work and patient analysis, but we believe that it is fair to say that 
these efforts failed to uncover any well-understood, predictable pattern.

In this paper and a subsequent part we shall describe in such cases 
the GIT-compactification \emph{as a stratified space} in terms of the 
BB-compactification plus some simple data of an arithmetic nature. 
Not only will this render the stratified structure more transparant, but
we will also see that the theory is strong enough to be able to guess in 
many cases what the GIT-compactification should be like.

We shall be more explicit now.
Let $\BB$ be a complex ball, $\LL\to\BB$ its natural automorphic line bundle
and $\G$ an arithmetic group operating properly on $(\BB,\LL)$. Then the orbit space $X:=\G\bs \BB$ is a quasiprojective variety and 
$\LL$ descends to an `orbiline bundle' $\Lcal$ over $X$. 
The Baily-Borel theory projectively compactifies $X$ by adding only finitely many points (the cusps).  
Suppose that we are also given a $\G$-invariant locally finite union 
of hyperballs in $\BB$. This determines a 
hypersurface in $X$ whose complement we shall denote by $X^\circ$. 
The central construction in this paper is a natural 
projective compactification $\widehat{X^\circ}$ of $X^\circ$
that (i) is explicitly given as a blowup of the Baily-Borel compactication
of $X$ followed by a blowdown and (ii) is naturally stratified.   
Perhaps the most important property for applications is that the boundary 
$\widehat{X^\circ}-X^\circ$ is everywhere of codimension $\ge 2$ in case
every nonempty intersection of the given hyperballs (including the intersection 
with empty index set, i.e., $\BB$) has dimension $\ge 2$.

Here is how this is used.
Let $G$ be a  reductive group acting on an integral 
quasi-projective variety endowed with an ample bundle 
$(Y,\eta)$. Assume that we are given a nonempty open subset
$U\subset Y$ which is a union of $G$-stable orbits. Then 
$G\bs (U, \eta |U)$ exists as a quasi-projective variety with orbiline 
bundle. To be precise, the $G$-invariants in the algebra of sections of powers of $\eta$ is finitely generated and the normalized proj of this 
graded algebra is 
a projective variety whose points correspond to the minimal orbits 
in the semistable part $Y^\ss$ (and is therefore denoted $G\bss Y^\ss$) and
contains $G\bs U$ as an open-dense subvariety.
Assume that we are given an isomorphism of orbiline bundles
\[
G\bs (U,\eta |U)\cong (X^\circ,\Lcal|X^\circ)
\]
for some triple $(X,X^\circ, \Lcal)$ as above. Assume also that the 
boundary at each side (so in $G\bss Y^\ss$ resp.\ $\widehat{X^\circ}$) is of codimension $\ge 2$. Perhaps contrary to what one might think, there are many examples  of interest for which these assumptions are satisfied with the isomorphism in question then being given by a period map.
Then we show that the isomorphism $G\bs U\cong X^\circ$ extends to an isomorphism
\[
G\bss Y^\ss\cong \widehat{X^\circ}.
\] 
It turns out that in all cases that we are aware of the stratification of 
$\widehat{X^\circ}$ has an interpretation in terms of the left hand side. 
Among these are the moduli space of quartic curves (equivalently, 
of degree two Del Pezzo surfaces) and of rational elliptic fibrations 
admitting a section (equivalently, of degree one Del Pezzo surfaces).\\ 

In the course of our discussion we touch briefly on a few topics that are 
perhaps not indispensable for getting at our main results, but that we 
included here since they involve little extra effort, 
can be illuminating to put things in perspective, and have an 
interest on their own.

One of these starts with the observation that the construction of 
$\widehat{X^\circ}$ makes sense and is useful in
a wider setting than ball quotients. For instance, we may take for 
$X^\circ$ the complement of an arrangement in a projective space or more
generally, in a torus embedding. This case is relevant for producing 
a compactification of the universal smooth genus three curve and the universal
smooth cubic surface, as they involve adjoint tori of type $E_7$ and $E_6$ 
respectively. We showed on previous occasions (\cite{looij1}, \cite{looij3})
that these constructions are sufficiently explicit for calculating the 
orbifold fundamental group of these universal objects. (This enabled us to find
a new, natural, simple presentation of the pointed mapping class group of 
genus three.)

Another such topic concerns a relatively simple nontrivial necessary 
condition that an `arithmetic' locally finite union of hyperballs in 
a complex ball must satisfy in order that it be the zero set of
an automorphic form (Section \ref{productexp}).\\

Let us finally point out that this paper builds on work of ours that 
goes back to the eighties, when we set out to understand the 
semi-universal deformations of triangle singularities and the GIT 
compactifications of J.~Shah of polarized K3 surfaces of degree 2 and 4.
Our results were announced in \cite{looij2}, but their technical 
nature was one of the reasons that the details got only partially 
published (in a preprint form in \cite{looij:toric} and in H.~Sterk's 
analysis of the moduli space of Enriques surfaces \cite{sterk}). We 
believe that the situation has now changed. Until our recent work 
with Gert Heckman \cite{hl} we had not realized that the constructions 
also work for ball quotients and that they form an attractive class to 
treat before embarking on the more involved case of type IV domain 
quotients. It also turned out that by doing this case first helped us 
in finding a natural setting for the construction (in particular, the notion \ref{def:linear} of a linearization of an arrangement thus presented itself) 
which makes its relative simplicity more visible. This was for us an important 
stimulus to return to these questions.
In part II we shall develop the story for type IV domains.

\medskip
I thank Gert Heckman for his many useful comments on the first version of this paper.

\section{Arrangements and their linearizations}

We begin with a definition.

\begin{definition}
Let $X$ be a connected complex-analytic manifold and $\Hcal$
a collection of smooth connected hypersurfaces of $X$.
We say that $(X,\Hcal)$ is an \emph{arrangement}
if each point of $X$ admits a coordinate neighborhood such that every $H
\in\Hcal$ meeting that neighborhood is given by a linear equation. 
In particular, the collection $\Hcal$ is locally finite.
\end{definition}

Fix an arrangement $(X,\Hcal)$. Suppose we are given
a line bundle $\Lcal$ on $X$ and for every $H\in\Hcal$ an isomorphism
$\Lcal (H)\otimes\Ocal_H\cong \Ocal_H$,
or equivalently, an isomorphism between the normal bundle 
$\nu_{H/X}$ and the coherent restriction of $\Lcal^*$ to $H$.
We shall want these isomorphisms to satisfy a
normal triviality condition which roughly says that 
if $L$ is any connected component of an intersection of members of 
$\Hcal$, then these isomorphisms trivialize the 
projectivized normal bundle of $L$. To be precise, if
$\Hcal^L$ denotes the collection of $H\in\Hcal$ that contain $L$, then the  
normal bundle $\nu_{L/X}$ is naturally realized as a subbundle of 
$\Lcal^*\otimes\Ocal_L^{\Hcal^L}$ by means of the homomorphism
\[
\nu_{L/X}\to \oplus_{H\in\Hcal^L} \nu_{H/X}\otimes\Ocal_{L}\cong
\Lcal^*\otimes\Ocal_L\otimes_\CC \CC^{\Hcal^L}.
\]
The condition alluded to is given in the following 

\begin{definition}\label{def:linear}
A \emph{linearization} of an arrangement 
$(X,\Hcal)$ consists of the data of a line bundle $\Lcal$ on $X$ 
and for every $H\in\Hcal$ an isomorphism $\Lcal (H)\otimes\Ocal_H\cong \Ocal_H$ 
such that the image of the above embedding 
$\nu_{L/X}\to \Lcal^*\otimes\Ocal_L\otimes_\CC \CC^{\Hcal^L}$
is given by a linear subspace $N(L,X)$ of $\CC^{\Hcal^L}$.
(Here $L$ is any connected component of an intersection of members of 
$\Hcal$.) We then refer to $N(L,X)$ as the \emph{normal space} of $L$ in $X$ 
and to its projectivization $\PP (L,X):=\PP (N(L,X))$ as the 
\emph{projectively normal space} of $L$ in $X$.
\end{definition}

The condition imposed over $L$ is 
empty when $L$ has codimension one and is automatically satisfied
when each $H$ is compact: if $I\subset\Hcal^L$ consists of $\codim L$   
elements such that $L$ is a connected component of $\cap_{H\in I} H$,   
then $\nu_{L/X}$ projects isomorphically on the subsum of 
$\Lcal^*\otimes \CC^{\Hcal^L}$ defined by $I$ and hence $\nu_{L/X}$
is given by a matrix all of whose entries are nowhere 
zero sections of $\Ocal_L$. So these entries are constant.

Here are some examples.

\begin{example}[Affine arrangements]\label{ex:affarr}
The most basic example is when $X$ is a complex affine space,
$\Hcal$ a locally finite collection of affine-linear hyperplanes and 
$\Lcal =\Ocal_X$. If the collection is finite, then by including the hyperplane 
at infinity, this becomes a special case of:
\end{example}
 
\begin{example}[Projective arrangements]\label{ex:projarr}
Here $X$ is a projective space $\PP (W)$,
$\Hcal$ is a finite collection of projective hyperplanes and 
$\Lcal =\Ocal_{\PP^n}(-1)$. 
\end{example}

\begin{example}[Toric arrangements]\label{ex:toric} 
Now $X$ is a principal homogeneous space of an algebraic torus $T$,  
$\Hcal$ is a finite collection of orbits of hypertori of $T$ in $X$ and
$\Lcal=\Ocal_X$. From this example we may obtain
new ones by extending these data to certain smooth
torus embeddings of $X$. The kernel of the exponential map
$\exp: \lie (T)\to T$ is a lattice $\lie (T)(\ZZ)$ in $\lie (T)$. 
So it defines a $\QQ$-structure on $\lie (T)$.
We recall that a normal torus embedding $X\subset X_\Sigma$ is given by 
a finite collection $\Sigma$ of closed rational polyhedral cones 
in $\lie (T)(\RR)\}$ with the property 
that the intersection of any two members is a common facet of these 
members. The torus $T$ acts on $X_\Sigma$ and the orbits of this action 
are naturally indexed by $\Sigma$: the orbit corresponding to
$\sigma\in\Sigma$ is identified with the quotient $T(\sigma)$ of $T$ by
the subtorus whose Lie algebra is the complex span of  
$\sigma$. The torus embedding is smooth precisely when each 
$\sigma$ is spanned by part of a basis of $\lie (T)(\ZZ)$ and in that case
$\Delta :=X_\Sigma -X$ is a normal crossing divisor.

If we want the closure of $H$ in $X_\Sigma$ 
to meet $\Delta$ transversally for every $H\in\Hcal$, then we must require 
that $\Sigma$ is closed under intersections with the real hyperplanes $\lie(T_H)(\RR)$, where $T_H\subset T$ is the stabilizer of $H$ (a hypertorus). Then the normal bundle of
$\Hbar$ is  isomorphic to $\Ocal (-\Delta )\otimes \Ocal_{H_\Sigma}$
(choose a general $T$-invariant vectorfield on $X_\Sigma$ and restrict 
to $\Hbar$). So $\Ocal (\Delta )$ may take the role of $\Lcal$.

A case of special interest is when $T$ is the adjoint torus of
a semisimple algebraic group $G$ and $\Hcal$ is the 
collection of hyperplanes defined by the roots. 
The corresponding hyperplanes in $\lie (T)(\RR)$
are reflection hyperplanes of the associated Weyl group and these decompose
the latter into chambers. Each chamber is spanned by 
a basis of $\lie (T)(\ZZ)$ and so the corresponding decomposition $\Sigma$
defines a smooth torus embedding $T\subset T_\Sigma$. Any root 
of $(G,T)$, regarded as a nontrivial character of $T$,  
extends to a morphism $T_\Sigma\to \PP^1$ that is smooth
over $\CC^\times$. The fiber over $1$ is the closure in $T_\Sigma$ 
of the kernel of this root and hence is smooth. 
\end{example}

\begin{example}[Abelian arrangements]\label{ex:abelian} 
Let $X$ be a torsor (i.e., a principal homogeneous space)
over an abelian variety, $\Hcal$ a collection of abelian 
subtorsors of codimension one and $\Lcal =\Ocal_X$.
\end{example}

\begin{example}[Diagonal arrangements]\label{ex:diagonal} 
Let be given a smooth curve $C$ of genus $g$ 
and a nonempty finite set $N$. 
For every two-element subset $I\subset N$  
the set of maps $N\to C$ that are constant on $I$ is a diagonal 
hypersurface in $C^N$ and the collection of these 
is an arrangement. But if $|N|>2$ and $C$ is a general complete 
connected curve of genus $\not= 1$, then a linearization will not exist:
it is straightforward to check that there is no linear combination 
of the diagonal hypersurfaces and pull-backs of divisors on the factors with 
the property that its restriction to every diagonal hypersurface is 
linearly equivalent to its conormal bundle. (If the genus is one, then 
the diagonal arrangement is abelian, hence linearizable.)
\end{example}

\begin{example}[Complex ball arrangements]\label{ex:ball}
Let $W$ be a complex vector space of finite dimension $n+1$
equipped with  a Hermitian form $\psi : W\times W\to \CC$ of 
signature $(1,n)$,
$\BB\subset \PP (W)$ the open subset defined $\psi (z,z)>0$ and 
$\Ocal_\BB (-1)$ the restriction of the `tautological' bundle 
$\Ocal_{\PP (W)}(-1)$ to $\BB$. Then $\BB$ is a complex ball of 
complex dimension $n$. 
It is also the Hermitian symmetric domain of the special unitary group 
$\SU (\psi)$. The group $\SU (\psi)$ acts on the line bundle 
$\Ocal_\BB (-1)$ and the latter is the natural automorphic line bundle 
over $\BB$. Every hyperplane $H\subset W$ of hyperbolic signature
gives a (nonempty) hyperplane section $\BB_H:=\PP(H)\cap \BB$ of $\BB$.
The latter is the symmetric domain of the special unitary group 
$\SU (H)$ of $H$. 
In terms intrinsic to $\BB$: $\BB_H$ is a totally geodesic
hypersurface of $\BB$ and any such hypersurface is of this form.
Notice that the normal bundle of $\BB_H$ is 
naturally isomorphic with $\Ocal_{\BB_H}(1)\otimes_\CC W/H$.
Hence we have an $\SU (H)$-equivariant isomorphism
\[
\Ocal_\BB (-1)(\BB_H)\otimes \Ocal_{\BB_H}\cong\Ocal_{\BB_H}\otimes_\CC W/H\cong
\Ocal_{\BB_H}.
\] 
Locally finite collections of hyperplane sections of $\BB$ arise 
naturally in an arithmetic setting.
Suppose that $(W,\psi)$ is defined over an imaginary 
quadratic extension $k$ of $\QQ$, which we think of as a subfield
of $\CC$. So we are given
a $k$-vector space $W(k)$ and an Hermitian form $\psi(k): W(k)\times 
W(k)\to k$ such that we get $(W,\psi)$ after extension of scalars.
We recall that a subgroup of $\SU (\psi )(k)$
is said to be \emph{arithmetic} if it is commensurable with the group 
of elements in $\SU (\psi)(k)$ that stabilize the $\Ocal_k$-submodule
spanned by a basis of $W(k)$ (where $\Ocal_k$ denotes the ring of integers
of $k$). It is known that an arithmetic subgroup 
of $\SU (\psi)(k)$ acts properly discontinuously on $\BB$.
Let us say that a collection $\Hcal$ of hyperbolic hyperplanes 
of $W$ is \emph{arithmetically defined} if $\Hcal$, regarded as a subset 
of $\PP (W^*)$, is a finite union of orbits of an arithmetic subgroup of 
$\SU (\psi)(k)$ and has each point defined over $k$. In that case  
the corresponding collection of hyperplane sections $\Hcal |\BB$ of 
$\BB$ is locally finite (this known fact is part of Lemma 
\ref{interiorgen} below).
\end{example}

\begin{example}[Type IV arrangements]\label{ex:ivarr}
Let $V$ be a complex vector space of dimension $n+2$ equipped with 
a nondegenerate symmetric bilinear form $\phi : V\times V\to \CC$.
Assume that $(V,\phi)$ is defined over $\RR$ in such a way that
$\phi$ has signature $(2,n)$.
Let $\DD\subset \PP (V)$ be a connected component of
the subset defined $\phi (z,z)=0$ and $\phi (z, \bar{z})>0$. We let
$\Ocal_\DD (-1)$ be the restriction of $\Ocal_{\PP (V)}$ to $\DD$ as in the 
previous example.
Then $\DD$ is a the Hermitian symmetric domain of the 
$\SL (\phi )(\RR)$-stabilizer of $\DD$. 
The group $G(V)$ also acts on the line bundle 
$\Ocal_\DD (-1)$ and the latter is the natural automorphic line bundle 
over $X$. Every hyperplane $H\subset V$ defined over $\RR$ 
of signature $(2,n-1)$
gives a (nonempty) hyperplane section $\DD_H :=\PP(H_\CC)\cap \DD$ of $\DD$.
This is a symmetric domain for the group its $\SL (\phi )(\RR)$-stabilizer. 
As in the previous
example, it is a totally geodesic hypersurface of $\DD$, any such 
hypersurface is of this form, and we have an equivariant isomorphism
$\Ocal_\HH (-1)(\DD_H)\cong \Ocal_{\DD_H}$.
Any collection of totally geodesic hypersurfaces of $X$ 
that is arithmetically defined in the sense of the example above
(with the number field $k$ replaced by $\QQ$) is locally finite on $\DD$.

A ball naturally embeds in a type IV domain: if
$(W,\psi)$ is as in the ball example, then
$V:=W\oplus \overline{W}$ has signature $(2,2n)$. A real structure
is given by stipulating that the interchange map is complex conjugation
and this yields a nondegenerate symmetric bilinear form $\phi$ defined over
$\RR$. For an appropriate choice of component $\DD$, we thus get an 
embedding $\BB\subset \DD$. If $(W,\psi)$ is defined over an imaginary 
quadratic field, then  $(W\oplus\overline{W},\phi)$ is defined over $\QQ$ 
and an arithmetic arrangement on $\DD$ restricts to one on $\BB$.
\end{example}
 
\section{Blowing up arrangements}

\subsection{Blowing up a fractional ideal}
In this subsection we briefly recall the basic notion of blowing up a 
fractional ideal. Our chief references are \cite{harts} and \cite{ega}. 

Suppose $X$ is a variety and $\Ical\subset \Ocal_X$ is 
a coherent ideal. Then $\oplus_{k=0}^\infty \Ical ^k$ is a 
$\Ocal_X$-algebra that is generated as such by its degree one summand (it should
be clear that $\Ical ^0:=\Ocal_X$). Its proj is a projective
scheme over $X$, $\pi :\tilde X\to X$, called the \emph{blowup of} $\Ical$, with the 
property that $\pi^*\Ical_X$ is invertible and relatively very ample.

If $X$ is normal and $\Ical$ defines a reduced subscheme of $X$, then its 
blowup is normal also. 
The variety underlying $\tilde X$ is over $X$ locally given as follows:
if $\Ical$ is generated over an open subset $U\subset X$ by 
$f_0,\dots ,f_r\in \CC [U]$, then these generators determine a rational map 
$[f_0:\cdots :f_r]: U\to\PP^r$ and the closure of the graph of this map in 
$U\times \PP^r$ with its projection onto $U$ can be identified with 
$(\tilde X_U)_\red\to U$. 
This construction only depends on $\Ical$ as a $\Ocal_X$-module. So if
$\Lcal$ is an invertible sheaf on $X$, and $\Ical$ is a nowhere zero 
coherent subsheaf of $\Lcal$, then we still have defined the blowup
of $\Ical$ as the proj of $\oplus_{k=0}^\infty \Ical ^{(k)}
\subset \oplus_{k=0}^\infty\Lcal^{\otimes k}$, where $\Ical ^{(k)}$ denotes
the $k$th \emph{power} of $\Ical$: the image of $\Ical^{\otimes k}$ in 
$\Lcal^{\otimes k}$. In fact, for the definition 
it suffices that $\Ical$ is a nowhere zero coherent subsheaf of the sheaf 
of rational sections of $\Lcal$. 

The coherent pull-back of $\Ical$ along $\tilde X\to X$ is invertible and the
latter morphism is universal for that property: any morphism from a scheme 
to $X$ for which coherent pull-back of $\Ical$ is invertible factorizes over $\tilde X$.
In particular, if  $Y\subset X$ is a closed subvariety, then 
the blowup $\tilde X\to X$ of the ideal defining $Y$ is an isomorphism
when $Y$ is a Cartier divisor. If $Y$ is the support  
of an effective  Cartier divisor, then $\tilde X\to X$ is still finite, but if 
$Y$ is only a hypersurface, then some fibers of $\tilde X\to X$ may have
positive dimension.

Here is a simple example that has some relevance to what will follow. The fractional
ideal in the quotient field of $\CC [z_1,z_2]$ generated by $z_1^{-1}z_2^{-1}$
is uninteresting as it defines the trivial blowup of $\CC^2$. But the blowup of the ideal generated by $z_1^{-1}$ and $z_2^{-1}$ amounts to the
usual blowup of the origin in $\CC^2$.

\subsection{The arrangement blowup} Let $(X,\Hcal)$ be an arrangement. 
There is a simple and straightforward way to find a modification 
$X_{\Hcal}\to X$ of $X$ such that the preimage of $D$ is a normal
crossing divisor: first blow up the union of the dimension  
$0$ intersections of members of the $\Hcal$, then  
then the strict transform of the dimension  
$1$ intersections of members of the $\Hcal$, and so on, finishing
with blowing up the strict transform of the dimension  
$n-2$ intersections of members of the $\Hcal$:
\[
X=X^0 \leftarrow X^1\leftarrow\cdots\leftarrow 
X^{n-2}=\tilde X^{\Hcal}.
\]
We refer to $\tilde X^{\Hcal}$ as the blow-up of $X$ defined by the arrangement
$\Hcal$ or briefly, as the \emph{$\Hcal$-blow-up} of $X$. 
To understand the full picture, it is perhaps best to do one blow-up at 
a time. In what follows we assume that we are also given a linearization 
$\Lcal$ of the arrangement $(X,\Hcal)$.
We begin with a basic lemma.

Let us denote by $\PO (\Hcal )$ the partially ordered set of connected 
components of intersections of members of $\Hcal$ 
which have positive codimension. For $L\in\PO(\Hcal )$, we denote  
by $\Hcal^L$ the collection of $H\in\Hcal$ containing $L$
as a \emph{lower} dimensional subvariety (so this collection 
is empty when $\codim (L)=1$). The following lemma is clear.

\begin{lemma}\label{normals} 
Given $L\in\PO (\Hcal )$, then the members of $\Hcal^L$ define a linear
arrangement in $N(L,X)$. If $L'\in \PO (\Hcal )$ contains $L$, then we
have a natural exact sequence of vector spaces 
\[
0\to N(L,L')\to N(L,X)\to N(L',X)\to 0.
\]
\end{lemma}

If $\Lcal$ is an invertible sheaf on $X$, then we shall write 
$\Lcal (\Hcal)$ for the subsheaf $\sum_{H\in\Hcal} \Lcal(H)$
of the sheaf of rational sections of $\Lcal$.
This sheaf is a coherent $\Ocal_X$-module, but 
need not be invertible. In particular, it should not be confused with 
the invertible sheaf $\Lcal(\sum_{H\in\Hcal} H)$. The $k$th power 
of $\Lcal (\Hcal)$ `as a fractional ideal' is
\[
\Lcal (\Hcal)^{(k)}:= \sum_{(H_1,\dots ,H_k)\in\Hcal^k} 
\Lcal^{\otimes k} (H_1+\cdots +H_k) 
\]
and so its blowup is the proj of $\oplus_{k=0}^\infty \Lcal (\Hcal)^{(k)}$.
\\

Let $L\subset X$ be a minimal member of $\PO (\Hcal )$ and of codimension 
$\ge 2$ and let $\pi : \tilde X\to X$ blow up $L$.
The exceptional divisor $E:=\pi^{-1}L$ is then the projectivized normal 
bundle $\PP (\nu_{L/X})$ of $L$. In case a linearization exists, the previous lemma identifies this exceptional divisor with $L\times \PP(L,X)$ in such a manner that the strict transform $\tilde H$ of $H\in\Hcal^L$ meets it in 
$L\times\PP(L,H)$.

\begin{lemma}\label{basicblowup}
Suppose that the arrangement $(X,\Hcal)$ is linearized by the invertible
sheaf $\Lcal$. Then the strict transform of $\Hcal$ in 
$\tilde X$ is an arrangement that is naturally linearized by
$\tilde\Lcal :=\pi^*\Lcal (E)$. Precisely, if $\boxtimes$
stands for the exterior coherent tensor product, then
\begin{enumerate}  
\item[(i)] $\nu_{E/\tilde X}=
(\Lcal^*\otimes \Ocal_L)\boxtimes \Ocal_{\PP(L,X)}(-1)$, 
\item[(ii)] 
$\pi^*\Lcal (\Hcal)\otimes\Ocal_E=
(\Lcal^*\otimes \Ocal_L)\boxtimes 
(\sum_{H\in \Hcal^L}\Ocal_{\PP(L,X)}(-1)(\PP(L,H))$ and 
\item[(iii)] for every $H\in\Hcal^L$ we have 
$\tilde\Lcal (\tilde H)\otimes \Ocal_{\tilde H}\cong 
\Ocal_{\tilde H}$.
\end{enumerate}
\end{lemma}
\begin{proof} 
Assertion (i) is an immediate consequence of Lemma \ref{normals}.
If $H\in\Hcal^L$, then the coherent restriction of 
$\tilde\Lcal (\tilde H)=(\pi^*\Lcal)(E+\tilde H)=\pi^*(\Lcal (H))$  
to the preimage of $H$ is trivial. Hence the same is true for its
coherent restriction to $E$ or $\tilde H$, which proves assertion (ii).
It is clear that $\Ocal_{\tilde X}(\tilde H )\otimes\Ocal_E=
\Ocal_L\boxtimes\Ocal_{\PP(L,X)}(\PP(L,H))$. Now (iii) follows also: 
\begin{align*}
\pi^*\Lcal (\Hcal)\otimes\Ocal_E
&=(\pi^*\Lcal )(\sum_{H\in\Hcal^L} \Ocal_E(E+\tilde H))\\
&=\pi^*\Lcal \otimes\nu_{E/\tilde X}(\sum_{H\in\Hcal^L}\Ocal_E(E\cap\tilde H)\\
&=(\Lcal^*\otimes \Ocal_L)\boxtimes 
(\sum_{H\in\Hcal^L} \Ocal_{\PP(L,X)}(-1)(\PP(L,H)).
\end{align*}
\end{proof}

Hence  $(\tilde X,\{ \tilde H\}_{H\in\Hcal}, \tilde\Lcal)$ satisfies
the same hypotheses as $(X,\{ H\}_{H\in\Hcal},\Lcal)$.
Our gain is that we have eliminated an intersection component of 
$\Hcal$ of minimal dimension. We continue in this manner,
until the strict transforms of the members of $\Hcal$ are disjoint
and end up with $\tilde X^\Hcal\to X$.

\begin{convention}\label{convention}
If an arrangement $\Hcal$ on a connected complex manifold 
$X$ is understood, we often omit it from notation that a priori
depends on $\Hcal$. For instance, we may write $\tilde X$ for the 
corresponding blow-up $\tilde X^\Hcal$ of $X$ and $X^\circ$ for the 
complementary Zariski open subset (as lying either in $X$ or in $\tilde X$). 
\end{convention}

The members of $\Hcal-\Hcal^L$ that meet $L$ define an 
arrangement $\Hcal_L$ in $L$. So we have defined $\tilde L$
and $L^\circ$. It is realized as the strict transform of $L$ 
under the blowups of members of $\PO (\Hcal )$ smaller than $L$.
Lemma \ref{basicblowup} shows that the projectivized normal bundle of 
$\tilde L$ in $X$ can be identified with the trivial bundle 
$\tilde L\times \PP(L,X)$ such that the members 
of $\Hcal^L$ determine a projective arrangement in $\PP(L,X)$. 
So we have defined $\tilde{\PP}(L,X)$ and $\PP(L,X)^\circ$.
The preimage of $\tilde L$ in $\tilde X$ is a smooth divisor 
that can be identified with 
\[
E(L):=\tilde L\times \tilde \PP(L,X).
\]
It contains
\[
E(L)^\circ:=L^\circ\times \PP(L,X)^\circ.
\]
as an open-dense subset. 

Lemma \ref{basicblowup} yields with induction:

\begin{lemma}\label{basicidealblowup}
Under the assumtions of Lemma \ref{basicblowup} we have:
\begin{enumerate}
\item[(i)] The morphism $\pi :\tilde X\to X$ is obtained by 
blowing up the fractional ideal sheaf $\Ocal (\Hcal)$, 
\item[(ii)] for every $L\in \PO (\Hcal )$, 
the coherent pull-back of $\Lcal (\Hcal)$ to $E(L)$ 
is identified with the coherent pull-back of the sheaf 
$\Ocal_{\PP(L,X)}(-1)(\Hcal^L)$ 
on $\PP(L,X)$ under the projection 
$E(L)\to \tilde \PP(L,X)\to \PP(L,X)$ and
\item[(iii)] the coherent pull-back of $\Lcal (\Hcal)$ to $\tilde L$ 
is identified with the coherent pull-back of 
$\Lcal (\Hcal_L)$ to $\tilde L$ (a line bundle by (i)) tensorized over $\CC$
with the vector space $H^0(\Ocal_{\PP(L,X)}(-1)(\Hcal^L))$.
\end{enumerate}
\end{lemma}

\begin{lemma}\label{incidence}
If $L$ and $L'$ are distinct members of $\PO (\Hcal )$, 
then $E(L)$ and $E(L')$ do not intersect unless $L$ and $L'$ 
are incident. If they are, and $L\subset L'$, say, then 
\begin{align*}
E(L)\cap E(L') &=\tilde L\times \tilde{\PP}(L,L')\times\tilde {\PP}(L',X).
\end{align*}
\end{lemma}
\begin{proof}
The first assertion is clear.
The strict transform of $L$ in the blowing ups of members of $\PO (\Hcal )$ 
contained in $L$ is $\tilde L\times \PP(L,X)$. The strict transform
of $L'$ meets the latter in $\tilde L\times \PP(L,L')$. Then
$E(L)\cap E(L')$ must be the closure of the preimage of 
$\tilde L\times \PP(L,L')^\circ$ in $\tilde X$. If we next 
blow up the members of $\PO (\Hcal )$ strictly between $L$ and $L'$,
the strict transform of $\tilde L\times \PP(L,L') $ is 
$\tilde L\times \tilde{\PP}(L,L')$. Blowing up $L'$ yields 
$\tilde L\times \tilde{\PP}(L,L')\times \PP (L',X)$, Finally,
blowing up the members of $\PO (\Hcal )$ strictly containing $L'$ gives
$\tilde L\times \tilde{\PP}(L,L')\times \tilde{\PP}(L',X)$.
The lemma now follows easily.
\end{proof} 

This generalizes in a straightforward manner as follows:

\begin{lemma}\label{stratification}
Let $L_0,L_1,\dots ,L_r$ be distinct members of 
$\PO (\Hcal )$, ordered by dimension:
$\dim (L_0)\le \dim (L_1)\le\cdots\le \dim (L_r)$. Then
$\cap_i E(L_i)$ is nonempty if and only if these members make up a 
flag  $L_\pt=(L_0\subset L_1\subset\cdots\subset L_r)$ and in 
that case the intersection in question equals
\[
E(L_\pt):= \tilde L_0\times \tilde{\PP}(L_0,L_1)\times\cdots
\times \tilde{\PP}(L_{r-1},L_r)\times \tilde{\PP}(L_r,X).
\]
Moreover, its open-dense subset 
\[
E(L_\pt)^\circ := L_0^\circ\times \PP (L_0,L_1)^\circ\times\cdots
\times \PP (L_{r-1},L_r)^\circ\times \PP (L_r,X)^\circ
\]
is a minimal member of the Boolean algebra of subsets of 
$X_\Hcal$ generated by the divisors $E(L)$. These elements define a 
stratification of $\tilde X-X^\circ$.
\end{lemma}

\subsection{Minimal (Fulton-MacPherson) blowup of an arrangement}\label{minimalblowup}
The simple blowing-up procedure described above turns the arrangement into a normal crossing divisor, but
it is clearly not minimal with respect that property. If $V$ is
a vector space and $(W_i\subset V)_i$ is a collection subspaces
that is linearly independent in the sense that $V\to \oplus_i V/W_i$
is surjective, then the blowup of
$\cup_i W_i$ defines a morphism $\tilde V\to V$ for which the 
total transform of $\cup_i W_i$ is a normal crossing divisor:
the strict transform $\tilde W_i$ of $W_i$ is a smooth divisor 
and these divisors meet transversally, their common intersection
being fibered over $\cap_i W_i$ with fiber $\prod_i \PP (V/W_i)$.
It is also obtained by blowing up the strict transforms of the 
$W_i$'s in any order. 
So in case of an arrangement we can omit the blowups along 
$L\in\PO (\Hcal)$ with the property that the minimal 
elements of $\Hcal^L$ are independent along $L$ in the above sense.
Since we still end up with a normal crossing situation, we call 
this the \emph{minimal normal crossing resolution} of the arrangement. 
It is not hard to specify a fractional ideal on $X$ whose blowup yields
this minimal normal crossing resolution. These and similar blowups have
been introduced and studied by Yi Hu \cite{hu} who built on earlier work
of A.~Ulyanov \cite{ulyanov}.

A case that is perhaps familiar is the diagonal arrangement 
of Example \ref{ex:diagonal}: we get the Fulton-MacPherson  
compactification of the space of injective maps from a given
finite nonempty set to a given smooth curve $C$. (We may 
subsequently pass to 
the orbit space relative to the action automorphism group of $C$   
if that group acts properly.)

\section{Contraction of blown-up arrangements}
Throughout the rest of this paper, $(X,\Hcal,\Lcal)$ stands for a 
linearized arrangement.

\medskip
According to Lemma \ref{basicidealblowup} the pull-back 
$\pi^*\Lcal (\Hcal)$ on $\tilde X$ is 
invertible. Suppose for a moment that $X$ is compact and 
$\Lcal (\Hcal)$ is generated by its 
sections. Then so are $\tilde X$ and $\pi^*\Lcal (\Hcal)$ and hence the 
latter defines a morphism from $\tilde X$ to a projective space. 
By Lemma \ref{basicidealblowup} again, this morphism will be constant on the
fibers of the projections $E(L)\to \tilde{\PP}(L,X)$. In other words, 
if $R$ is the equivalence relation on $\tilde X$ generated by these
projections, then the morphism will factorize through the 
quotient space $\tilde X/R$. Our goal is to find conditions under which 
$\pi^*\Lcal (\Hcal)$ is semiample, more precisely, conditions
under which this quotient space is projective and $\pi^*\Lcal (\Hcal)$
is the coherent pull-back of an ample line bundle on $\tilde X/R$.

Let us first focus on $R$ as an equivalence relation. We begin with noting
that on a stratum $E(L_0\subset\cdots\subset L_r)^\circ$ the 
equivalence relation is defined by the projection on the last factor
$\PP (L_r,X)^\circ$. 
So the  quotient $\tilde X/R$ is as a set the disjoint union of 
$X^\circ$ and the projective arrangement complements $\PP (L,X)^\circ$, 
$L\in\PO (\Hcal )$. 

\begin{lemma}
The equivalence relation $R$, when viewed as a subset of 
$\tilde X\times\tilde X$, is equal to the union of the diagonal and the 
$E(L)\times_{\tilde{\PP}(L,X)} E(L)$. In particular, $R$ is closed
and hence $\tilde X/R$ is Hausdorff.
\end{lemma}
\begin{proof}
By definition $R$ contains the union in the statement. The opposite inclusion
also holds: $x\in E(L_0\subset\cdots\subset L_r)^\circ$ and 
$x'\in E(L'_0\subset\cdots\subset L'_s)^\circ$ are related if and only
if $L_r=L'_s$ and both have the same image in $\PP (L_r,X)^\circ$. 
\end{proof} 

\begin{theorem}\label{basicmodification}
Suppose $X$ compact and that some positive power
of $\Lcal (\Hcal)$ is generated by its sections 
and that these sections separate the points of $X^\circ$. Then the pull-back of
$\Lcal (\Hcal)$ to $\tilde X$ is a semi-ample invertible sheaf and a 
positive power of this sheaf defines a morphism from $\tilde X$ to a projective
space whose image $\hat X$ realizes the quotient space $X/R$ with the  
strata $X^\circ$ and $\PP(L,X)^\circ$ of $X/R$ realized as subvarieties 
of $\hat X$.
\end{theorem} 
\begin{proof} Let $\Lcal(\Hcal)^{(k)}$ be generated by 
its sections and separate the points of $X^\circ$. Then 
the invertible sheaf $\pi^*(\Lcal(\Hcal)^{(k)})$ has the same property. 
Let $f:\tilde X\to \PP^N$
be the morphism defined by a basis of its sections.
A possibly higher power of $\Lcal (\Hcal)$ will have a 
normal variety as its image and so it suffices to prove that $f$ 
separates the points of every stratum  of $\tilde X/R$.
By hypothesis this is the case for $X^\circ$. 
Given $L\in\PO (\Hcal )$ of codimension
$\ge 2$, then it follows from Lemma \ref{basicidealblowup} that the pull-back of 
$\Lcal (\Hcal)$ to $\tilde L$ is isomorphic to a line bundle
on $\tilde L$ tensorized over $\CC$ with 
$H^0(\Ocal_{\PP(L,X)}(-1)(\Hcal^L))$.
Any set of generating sections of $\Lcal (\Hcal)$ must therefore 
determine a set of generators of the vector space 
$H^0(\Ocal_{\PP(L,X)}(-1)(\Hcal^L))$. Such a generating set separates 
the points of $\PP(L,X)^\circ$. The same is true for a set of 
generating sections of a power of $\Lcal (\Hcal)$.
So for some positive $k$, $\Lcal(\Hcal)^{(k)}$ separates the points of every
stratum of $\tilde X/R$.
\end{proof}

This theorem applies to some of the examples mentioned at the beginning.

\subsection{Projective arrangements}
The hypotheses of Theorem \ref{basicmodification} 
are satisfied if the collection 
$\Hcal$ of hyperplanes has no point in common: if $f_H\in V^*$
is a defining equation for $H\subset\PP(V)$, then 
$\Ocal_{\PP(V)}(-1)(\Hcal)$ is generated by the sections 
$\{f_H^{-1}\}_{H\in\Hcal}$. There is no need to pass
to a higher power of this sheaf: the  variety 
$\hat{\PP}(V)$ appears as the image of the rational map
$\PP(V)\to \PP (\CC^\Hcal )$ with coordinates 
$(f_H^{-1})_{H\in\Hcal}$. If $V=\CC^{n+1}$ and 
$\Hcal$ is the set of coordinate hyperplanes, then 
$\hat{\PP^n}=\PP^n$ and this map is just the standard Cremona 
transformation in dimension $n$.

Of particular interest are the cases when $\Hcal$ is the set of 
reflection hyperplanes of a finite complex reflection group $G$ with
$V^G=\{ 0\}$. Then $\hat{\PP}(V)$ has a $G$-action so that we can form
the orbit space $G\bs\hat{\PP}(V)$. 
(A theorem of  Chevalley implies that $G\bs {\PP}(V)$
is rational, hence so is this completion.) 

\subsection{Toric arrangements defined by root systems}
Let $T$ be the adjoint torus of a semisimple algebraic group.
It comes with an action of the Weyl group $W$. 
The collection of reflection hyperplanes of $W$ in 
$\lie (T)(\RR)$ defines a torus embedding $T\subset T_\Sigma$. The roots define 
a collection of `hypertorus embeddings' $\Hcal$ in $T_\Sigma$. 
We may also state this as follows: the set $R$ of roots define an embedding
$T\to (\CC^\times)^R$, $T_\Sigma$ is the closure of its image in
$(\PP^1)^R$ and $H_\alpha$ is the fiber over $[1:1]$ of the 
projection of $T_\Sigma\to \PP^1_\alpha$ on the 
factor indexed by $\alpha$.
A $T$-orbit $T(\sigma)$ in $T\subset T_\Sigma$ is 
also an adjoint torus (associated to a parabolic subgroup), its closure
in $T_\Sigma$ is the associated torus embedding (hence smooth)
and the roots of $T$ restrict to roots of $T(\sigma)$ and vice versa. Denote by  $\Pi$ the collection of indivisible elements of $\lie (T)(\ZZ)$
that generate a (one dimensional) face of $\Sigma$ (these are just the 
coweights that are fundamental relative a chamber).
The closure $D(\varpi)$ of $T(\RR_{>0}\varpi)$ in $T_\Sigma$ is 
a (smooth) $T$-invariant divisor and all such are so obtained.
The meromorphic functions $\{ f_\alpha:=(e^\alpha -1)^{-1}\}_{\alpha\in R}$ 
separate the points
of $T^\circ$. They also separate the $T$-orbits from each other: 
The divisor of $f_\alpha$ is
\[
-H_\alpha+
\sum_{\sigma\in\Sigma_1} \alpha (\varpi)D(\varpi).
\]
Two distinct members of $\Pi$ can be distinguished by a root and so the corresponding $T$-orbits are separated 
by the corresponding meromorphic function. Given a $\varpi\in\Pi$, 
then each $f_\alpha$ with $\alpha (\varpi)=0$ restricts to a 
meromorphic function on $D(\varpi)$ of the same type. For that reason 
these meromorphic functions also separate the points of 
any stratum $T(\sigma)^\circ$.

So $\widehat{T_\Sigma}$ exists as a projective variety. It comes with
an action of the Weyl group $W$ of the root system.
Notice that the subvariety of $\widehat{T_\Sigma}$ over the 
identity element of $T$ is the modification $\hat{\PP}(\lie (T))$ 
of the projectivized tangent space 
$\PP (\lie (T))$ defined by the arrangement of tangent spaces of 
root kernels.

\begin{example}[The universal cubic surface]
Suppose $(T,W)$ is an adjoint torus  of type $E_6$.
Then the construction of \cite{looij1} (see also \cite{looij3}) shows 
that $(\{\pm 1\}\cdot W)\bs\widehat{T_\Sigma}$ may be 
regarded as a compactification of the universal smooth cubic surface. Let us
briefly recall the construction. If $S\subset \PP^3$ is a cubic surface, then
for a general $p\in S$, the projective tangent space of $S$ at $p$ meets 
$S$ in a nodal cubic $C_p$. The latter is also an anticanonical curve.
The identity component $\pic_o(C_p)$ of the  Picard group of $C_p$ 
is a copy of $\CC^\times$. The preimage $\pic_o(S)$ of $\pic_o(C_p)$ in
$\pic (S)$ is a root lattice of type $E_6$ and so 
$\Hom (\pic_o(S),\pic_o(C_p))$ is an adjoint $E_6$-torus. The 
restriction map defines an element of this torus. It turns out that it does not
lie in any reflection hypertorus. We proved in
\cite{looij1} that this element is a complete invariant of the pair $(S,p)$: the $\{\pm 1\}\cdot W$-orbit in $T^\circ$ defined by the restriction map
determines $(S,p)$ up to isomorphism. All orbits in $T^\circ$ so arise and so 
$(\{\pm 1\}\cdot W)\bs T^\circ$ is a coarse moduli space for pairs $(S,p)$ as above.  

If we allow $p$ to be an arbitrary point of a smooth cubic, then $C_p$ 
may degenerate into any reduced cubic curve. 
The type of this curve corresponds in fact to  
$\{\pm 1\}\cdot W$-invariant union of strata of $\widehat{T_\Sigma}$
in a way that $(\{\pm 1\}\cdot W)\bs T^\circ$ is going to contain the coarse
moduli space of pairs $(S,p)$ with $p$ an arbitrary point of the smooth 
cubic $S$. To be precise: if $C_p$ becomes a cuspidal curve (with a cusp at $p$), 
then there is a unique stratum, it is projective of type $E_6$  (it is the one defined by the identity element of $T$: a copy of $\PP(\lie (T))^\circ$). If $C_p$ becomes a conic plus a line, then
we get the projective or toric strata of type $D_5$, 
depending on whether or not the line is tangent to the conic. Finally, if
$C_p$ becomes a sum of three distinct lines, then we get 
projective or toric strata of type $D_4$, depending on whether or not
the lines are concurrent.  
\end{example}

\begin{example}[The universal quartic curve]
The story is quite similar to the preceding case (we refer again to 
\cite{looij1} and \cite{looij3}). 
We now assume that $(T,W)$ is an adjoint torus of type $E_7$.
A smooth quartic curve $Q$ in $\PP^2$ determines a Del Pezzo surface of
degree $2$: the double cover $S\to\PP^2$ ramified along $Q$. This sets up a bijective correspondence between isomorphism classes of quartic curves and
isomorphism classes of Del Pezzo surface of degree $2$. If $p\in Q$
is such that the projective line $T_pQ$ meets $Q$ in two other distinct 
points, then the preimage $Q_p$ of $T_pQ$ in $S$ is a nodal genus one curve.
Starting with the homomorphism $\pic (S)\to\pic (Q_p)$, we define
$\pic_o(S)$ as before (it is a root lattice of type $E_7$) so that
we have an adjoint $E_7$-torus $\Hom (\pic_o(S),\pic_o(Q_p))$. 
Proceeding as in the previous case we find that 
$W\bs T^\circ$ is a coarse moduli space for pairs $(Q,p)$ as above
(now $-1\in W$) and that allowing $p$ to be arbitrary yields an open-dense 
embedding of the corresponding coarse moduli space in     
$W\bs\widehat{T_\Sigma}$. The added strata are as follows: allowing
$T_pQ$ to become a flex (but not a hyperflex) point yields the projective 
stratum of type $E_7$ over the identity element. If $T_pQ$ is a bitangent, 
then we get the toric or projective strata of type $E_6$, 
depending on whether the bitangent is a hyperflex.
Our compactification of $W\bs T^\circ$ has another interesting
feature as well. 
We have also $A_6$ and $A_7$-strata, both of projective type. These are single orbits and have an interpretation as the coarse moduli space of pairs $(Q,p)$,
where $p$ is a point of a hyperelliptic curve $Q$ of genus $3$: we are on a
$A_6$ or $A_7$-stratum depending on whether or not $p$ is a Weierstra\ss\ point.
So $W\bs\widehat{T_\Sigma}$ contains in fact the coarse moduli space of
smooth pointed genus three curves $\Mcal_{3,1}$. We used aspects of this 
construction in \cite{looij2} and \cite{looij3} to compute the rational 
cohomology resp.\ the orbifold fundamental group of $\Mcal_{3,1}$.
\end{example}
 
\subsection{Abelian arrangements}\label{abarr} Let $X$ be a torsor 
over an abelian variety $A$, 
and $\Hcal$ a finite collection of abelian subtorsors of codimension one.
Denote by $A_0$ the identity component of the group of 
translations of $X$ that stabilize $\Hcal$ (an abelian variety).
First assume that $A_0$ reduces to the identity element.
This ensures that $\Ocal_X(\sum_{H\in\Hcal} H)$ is ample. 
Hence a power of $\Ocal_X(\Hcal)$ separates the points of 
$X^\circ$. It is now easy to see that $\hat X$ is defined. 
For every $H\in\Hcal$, $X/H$ is an elliptic curve (the origin is 
the image of $H)$. A rational function on $X/H$ that 
is regular outside the origin and has a pole of order $k$ at the origin
can be regarded as a section of $\Ocal_X(k\Hcal)$. For a sufficiently
large $k$, the collection of these functions define a morphism from
$\tilde X$ to a projective space  which factorizes over 
$\hat X$. The morphism from $\hat X$ to this projective space is 
finite. In the general situation (where $A_0$ may have positive 
dimension), the preceding construction can be carried out in
a $A_0$-equivariant manner to produce a projective completion
$\hat X$ of $Y^\circ$ with $A_0$-action.

Here is a concrete example. 
If $R$ is a reduced root system with 
root lattice $Q$ and $E$ is an elliptic curve, then $X:=\Hom (Q,E)$
is an abelian variety on which the Weyl group $W$ of $R$ acts. The 
fixed point loci of reflections in $W$ define an abelian arrangement
$\Hcal$ on $X$ with $\cap_{H\in\Hcal} H$ finite. So  $\hat X$ is 
defined and comes with an action of $W$. 

\begin{example} 
For $R$ of type $E_6$,  
$(\{\pm 1\}\cdot W)\bs X^\circ$ is the moduli space of smooth cubic 
surfaces with hyperplane section isomorphic to $E$, provided that $E$ 
has no exceptional automorphisms. Hence $(\{\pm 1\}\cdot W)\bs \hat{X}$ 
is some compactification of this moduli space. 
\end{example}
\begin{example}
For $R$ of type $E_7$ 
there is a similar relationship with the moduli space of smooth quartic 
curves with a line section for which the four intersection 
points define a curve isomorphic to $E$.
\end{example}

\subsection{Ample line bundles over abelian arrangements}
This is a variation on \ref{abarr}, which we mention here
because of a later application to automorphic forms on ball quotients and
type IV domains with product expansions. 
Suppose that in the situation of \ref{abarr}
we are given an ample line bundle $\ell$ over $X$. Let 
$C$ be the corresponding affine cone over $X$, that is, the variety obtained from
the total space of the dual of $\ell$ by collapsing its zero section.
In more algebraic terms, $C_H$ is spec of the graded algebra
$\oplus_{k=0}^\infty H^0(\ell^{\otimes k})$. Each $H\in\Hcal$ defines
a hypersurface (a subcone of codimension one) $C_H\subset C$. We put
$C_\Hcal:=\cup_{H\in\Hcal} C_H$. 

\begin{lemma}\label{product} 
If $\Hcal\not=\emptyset$, then the hypersurface $C_\Hcal$ is 
the zero set of a regular function on $C$ if and only if 
the class of $\ell$ in $\pic (X)\otimes\QQ$ is a positive linear 
combination of the classes $H$, $H\in\Hcal$.
\end{lemma}
\begin{proof}
If $f$ is a regular function on $C$ with zero set $C_\Hcal$, then
$f$ must be homogeneous, say of degree $k$. So $f$ then defines
a section of $\ell$ with divisor $\sum_{H\in\Hcal} n_H H$, with
$n_H$ a positive integer. This  means that $\ell^{\otimes k}\cong
\Ocal_X(\sum_H n_H H)$. Conversely, an identity of this type implies the
existence of such an $f$.
\end{proof}

Since $\ell$ is ample, the condition that $\ell$ is a positive linear 
combination of the classes indexed by $\Hcal$, is not fulfilled if
$A_0$ (the identity component of the group of 
translations of $X$ that stabilize $\Hcal$) is positive dimensional.

\subsection{Dropping the condition of smoothness}
If in the situation of the previous theorem, $Y$ is a normal 
subvariety of $X$ which meets the members of $\PO (\Hcal )$ transversally, 
the strict transform
of $Y$ in $\hat X$ only depends on the restriction of the arrangement
to $Y$. This suggests that the smoothness condition imposed on $X$ 
can be weakened to normality. This is indeed the case: 

\begin{definition}
Let $X$ be a normal connected (hence irreducible) variety, 
$\Lcal$ be a line bundle over $X$. 
A \emph{linearized arrangement on $X$} consists of a locally finite 
collection of (reduced) Cartier divisors $\{ H\}_{H\in\Hcal}$ on $X$,
a line bundle $\Lcal$ and  
and for each $H\in\Hcal$ an isomorphism $\Lcal (H)\otimes\Ocal_H\cong\Ocal_H$  
such that for every connected component $L$ of an intersection of members
of $\Hcal$ the following conditions are satisfied:
\begin{enumerate}
\item[(i)] the natural homomorphism 
$\oplus_{H\in\Hcal^L}\Ical_H/ \Ical_H^2\otimes\Ocal_L\to \Ical_L/\Ical_L^2$
is surjective and not identically zero on any summand,
\item[(ii)] if we identify $\oplus_{H\in\Hcal^L}\Ical_H/ \Ical_H^2\otimes\Ocal_L$ with 
$\Lcal\otimes\Ocal_L\otimes_\CC \CC^{\Hcal^L}$ via the given isomorphisms,
then the kernel of the above homomorphism is spanned by a subspace 
$K(L,X)\subset \CC^{\Hcal^L}$ whose codimension is that of $L$ in $X$.
\end{enumerate}
\end{definition}

Notice that these conditions imply that the conormal sheaf 
$\Ical_L/\Ical_L^2$ is a locally free $\Ocal_L$-module of rank equal to 
the codimension of $L$ in $X$. 
The normal space $N(L,X)$ is now defined as the space of linear forms on 
$\CC^{\Hcal^L}$ that vanish on $K(L,X)$. 

The conditions of this definition are of course chosen in such a manner 
as to ensure that the validity of the discussion for the smooth case is not 
affected. In particular, we have defined
a blow-up $\tilde X\to X$ whose exceptional locus is a union of divisors
$E(L)\cong \tilde L\times\tilde \PP (L,X)$. 

\begin{theorem}\label{singversion} 
Suppose that is $X$ compact and that a  power of 
$\Lcal (\Hcal)$ is generated by its sections and separates the points of 
$X^\circ$. Then the pull-back of $\Lcal (\Hcal)$ to $\tilde X$ 
is semi-ample and the corresponding projective contraction 
$\tilde X \to \hat X$ has the property that $\hat X -X^\circ$ 
is naturally stratified into subvarieties indexed by $\PO (\Hcal)$: 
to $L\in \PO (\Hcal)$ corresponds a copy of $\PP (L,X)^\circ$ (this 
assignment reverses the order relation).
\end{theorem}

\section{Complex ball arrangements}

We take up Example \ref{ex:ball}. We shall assume that $n\ge 2$.
It is clear that Theorem \ref{basicmodification} does not
cover this situation and so there is work to do. 

\subsection{The ball and its natural automorphic bundle}
We choose an arithmetic subgroup $\G$ of $\SU (\psi)(k)$, although we
do not really want to fix it in the sense
that we always allow the passage to an arithmetic subgroup of $\G$ of finite 
index. In particular, we shall assume that $\Gamma$ is \emph{neat}, which means
that the multiplicative subgroup of $\CC^\times$ generated by the eigen values of 
elements of $\G$ is torsion free. This implies that
every subquotient of $\G$ that is `arithmetically defined'
is torsion free. 

We write $X$ for the $\G$-orbit space of $\BB$.
A \emph{cusp} of $\BB$ relative the given $k$-structure is an element 
in the boundary of $\BB$ in $\PP (W)$ that is defined over $k$. 
A cusp corresponds to an isotropic line $I\subset W$ defined over $k$.

Denote by $\LL^\times\subset W$ the set of $w\in W$ with
$\psi (w,w)>0$. The obvious projection $\LL^\times\to \BB$ is a
$\CC^\times$ bundle. We may view this as the complement of
the zero section of the tautological line bundle over $\PP (W)$ restricted
to $\BB$: 
\[
\LL :=\LL^{\times}\times ^{\CC^\times}\CC.
\] 
A a nonzero $v\in L$ defines a homogeneous function $f:L-\{ 0\}\to \CC$ of 
degree $-1$ characterized by the property that $f(v)=1$. So a nonzero holomorphic function on $\LL^\times$ that is homogeneous of degree $-k$ defines a holomorphic 
section of $\LL^{\otimes k}$ and vice versa.  

Let us point out here the relation between $\LL$ and the canonical bundle
of $\BB$. If $p\in\BB$ is represented by the line $L\subset W$, then
the tangent space $T_p\BB$ is naturally isomorphic to $\Hom (L,W/L)$.
So $\wedge^nT_p^*\BB\cong L^{n+1}\otimes\wedge^{n+1}W^*$. This proves
that the canonical bundle of $\BB$ is $\SU (\psi)$-equivariantly 
isomorphic to $\LL^{\otimes (n+1)}$. We regard $\LL$ as the natural 
automorphic bundle over $\BB$. Since $\G$ is neat, $\LL$ drops to 
a line bundle over $X$.

\subsection{The stabilizer of a cusp}\label{stabcusp}
The following bit of notation will be useful.

\begin{convention} 
Given a subspace $J$ of $W$ and a subset $\Omega\subset
W$ resp.\ $\Omega\subset \PP(W)-\PP(J)$, then 
$\Omega (J)$ denotes the image of $\Omega$ in $W/V$ resp.\
$\PP(W/J)$.
\end{convention} 

Suppose that $J$ is a degenerate subspace of $\BB$, so that
$I:=J\cap J^\perp$ is an isotropic line with $I^\perp =J+J^\perp$. 
Then
\[
\LL^\times (J)= W/J -I^\perp/J \quad \text{and hence}\quad 
\BB (J)=\PP(W/J)-\PP (I^\perp/J).
\]
In the first case we have the complement of a linear hyperplane in a 
vector space; in the second case the complement of a projective 
hyperplane in a projective space, hence an affine space
(with $\Hom (W/I^\perp ,I^\perp/J)$ as translation group). The former
is a $\CC^\times$-bundle over this affine space. Notice that for the 
maximal choice $J=I^\perp$, this is a $\CC^\times$-bundle over a 
singleton. 

Let $I\subset W$ be an isotropic line.
We begin with recalling the structure of the stabilizer 
$\SU (\psi)_I$. If $e$ is a generator of $I$,
then the unipotent radical $N_I$ of the $\SU (\psi)$-stabilizer of $I$ 
consists of transformations of the form
\[
T_{e,v}: z\in W\mapsto z+\psi (z,e)v-\psi (z,v)e
-\tfrac{1}{2}\psi (v,v)\psi (z,e)e 
\]
for some  $v\in I^\perp$. Notice that $T_{e,\lambda v}
=T_{\overline{\lambda} e,v}$ and that $T_{e,v+\lambda e}=T_{e,v}$ when
$\lambda\in\RR$. So $T_{e,v}$ only depends on the image of $e\otimes v$ 
in $\overline{I}\otimes_\CC I^\perp/(\overline{I}\otimes I)(\RR)$ 
(observe $\overline{I}\otimes I$ has a natural real structure indeed).
A simple calculation yields 
\[
T_{e,u}T_{e,v}= T_{e, u+v+\frac{1}{2}\psi (u,v)e}
\]
showing that $N_I$ is an Heisenberg group: its center
$Z(N_I)$ is identified with the real line $\sqrt{-1}(\overline{I}\otimes_\CC I)(\RR)$
and the quotient $N_I/Z(N_I)$ with the vector group $\overline{I}\otimes I^\perp/I$, 
the commutator given essentially by the skew form that the  
imaginary part $\psi$ induces on $I^\perp/I$. 
One can see this group act on $\BB$ as follows. 
The domain $\BB$ lies in $\PP(W)-\PP(I^\perp)$. The latter is an affine
space for $\Hom (W/I^\perp ,I^\perp)$ and so $N_I$ will act on this as an 
affine transformation group. This group preserves the fibration
\[
\PP(W)-\PP(I^\perp)\to \PP(W/I)-\PP(I^\perp/I)=\BB(I).
\]
This is a fibration by affine lines whose base $\PP(W/I)-\PP(I^\perp/I)$ 
is an affine space over $\Hom (W/I^\perp ,I^\perp/I)$.
The center $Z(N_I)$ of $N_I$ respects each fiber and acts there as 
a group of translations and the quotient vector group $N_I/Z(N_I)$ 
is faithfully realized on the base as its group of translations.
The domain $\BB$ is fibered over $\BB (I)$ by half 
planes such that the $Z(N_I)$-orbit space of a fiber is a half line. 

Let us first decribe $\BB$ as such in terms of (partial) coordinates:
choose $f\in V$ such that $\psi (e,f)=1$ and denote the orthogonal complement
of the span of $e$ and $f$ by $A$. It is clear that $A$ is negative definite
and can as an inner product space be identified with
$I^\perp/I$. The map $(s,a)\in \CC\times A\mapsto [se+f+a]\in \PP (V)-\PP (I^\perp)$ is an isomorphism of affine spaces and in these terms $\BB$ is defined by $2\text{Re}(s)>-\psi (a,a)$. This is known as the realization of $\BB$ as a Siegel domain of the second kind.

A somewhat more intrinsic description goes as follows.
Consider the group homomorphism $T_e:\CC\to\GL (W)$ defined by
\[
T^s_e(z)=z+s\psi (z,e)e.
\]
A simple computation shows that 
\[
\psi (T^s_ez, T^s_ez)=\psi (z,z)+2 \text{Re}(s)|\psi (z,e)|^2.
\]
The restriction of $T_e$ to the imaginary axis parametrizes $Z(N_I)$
and therefore $T_e$ maps to the complexification $\SU (\psi )$.
Since $T^s$ maps $\BB$ into itself for $\text{Re} (s)\ge 0$, the 
image $\Tcal_I$ in $GL(W)$ of the right half plane is a semigroup
that acts (freely) on $\BB$. As the notation suggests it does not depend
on the generator $e$ of $I$.
The fibers of $\BB\to\BB (I)$ are orbits of this semigroup.
The image of $T$ lies in the complexification of $\SU (\psi)$. It is 
normalized by $N_I$ and so $N_I\Tcal_I=\Tcal_IN_I$ is
a semigroup in $\SU (\psi)_\CC$. The domain $\BB$ is a free orbit of
this semigroup.    

Suppose now that $I$ is defined over $k$.
The Levi quotient $\SU (\psi)_I/N_I$ of $\SU (\psi)_I$
is the group of unitary
transformations of $I^\perp/I$, hence compact, and the corresponding 
subquotient of $\G$ is finite. Since we assumed $\G$ to be neat, 
the latter is trivial and so the stabilizer $\Gamma_I$ is contained in 
$N_I$. It is in fact a 
cocompact subgroup of $N_I$: the center $Z(\G_I)$ 
of $\G_I$ is infinite cyclic (a subgroup of 
$\sqrt{-1}(\overline{I}\otimes I)(\RR)$), and $\Gamma_I/Z(\G_I)$ 
is an abelian group that naturally lies as a lattice in the 
vector group $\overline{I}\otimes I^\perp/I$.

\subsection{Compactications of a ball quotient}\label{analyticstructure} 
Let $J$ be a complex subspace of $W$ with radical $I$.
The group $N^J$ of $\g\in\G$ that act as the identity on 
$J^\perp\cong \overline{W/J}^*$ is a normal (Heisenberg) subgroup of $N_I$  
with quotient identifyable with the vector group $\overline{I}\otimes I^\perp/J)$. The latter may be regarded
as the vector space of translations of the affine space
$\BB(J)=\PP(W/J)-\PP(I^\perp/J)$. Let us assume that $J$ is defined over $k$.
Then the discrete counterparts of these statements hold for $\G$. 
In particular, $\G_I/\G^J$ can be identified with
a lattice in $\overline{I}\otimes J/I$.
So if we denote the orbit space of this lattice
acting on $\BB(J)$ by $X(J)$, then $X(J)$ is a 
principal homogenenous space of a complex torus that has
$\overline{I}\otimes I^\perp/J)$ as its universal cover. For 
an intermediate $k$-space $I\subset J'\subset J$ we have a natural projection $X(I)\to X(J)$ of abelian torsors. Notice that $X(I^\perp)$ is 
just a singleton.

\begin{definition}
An \emph{arithmetic system} (of degenerate subspaces) assigns 
in a $\G$-equivariant manner to every $k$-isotropic line $I$ a  
degenerate subspace $j(I)$ defined over $k$ with radical $I$. We call the 
two extremal cases $\bb (I):=I^\perp$ resp.\ $\tor (I):= I$ the 
\emph{Baily-Borel system} and the \emph{toroidal system} respectively.
\end{definition}

Such a system $j$ leads to a compactification $X^j$ of $X$ as follows.
Form the disjoint unions
\begin{align*}
(\LL^\times)^j:=\LL^\times\sqcup 
\coprod_{I \text{ $k$-isotropic}} \LL^\times(j(I)) \quad \text{and}\quad
\BB^j:=\BB\sqcup 
\coprod_{I \text{ $k$-isotropic}} \BB (j(I)).
\end{align*}
Both come with an obvious action of $\G$ and the projection
$(\LL^\times)^j\to \BB^j$ is a $\CC^\times$-bundle.
We introduce a $\G$-invariant topology on these sets as follows.
Recall that $N^{j(I)}$ is the subgroup of $g\in N_I$
that act as the identity on $W/j(I)$. Then $N^{j(I)}\Tcal_I=
\Tcal_IN^{j(I)}$ is a semigroup and the fibers of
$\BB\to \BB(j(I))$ are orbits of this semigroup.
A basis of a topology on
$(\LL^\times)^j$ is specified by the following collection of subsets: 
\begin{enumerate}
\item[(i)] the open subsets of $\LL^\times$,
\item[(ii)] for every $k$-isotropic line $I$ and
$N^{j(I)}\Tcal_I$-invariant open subset $\Omega\subset \LL^\times$ 
the subset  $\Omega\sqcup \Omega(j(I))$.
\end{enumerate} 
The same definition with obvious modifications defines a topology for $\BB^j$. It makes
\[
\LL^j:=(\LL^\times)^j \times^{\CC^\times}\CC .
\]
a topological complex line bundle over $\BB^j$. We are interested in the
orbit space
\[
X^j:=\G\bs \BB^j
\]
and the line bundle over $X^j$ defined by $\G\bs \LL^j$. 
Notice that as a set, $X^j$ is the disjoint union of $X$ and the torsors
$X(j(I))$, where $I$ runs over a system of representatives of
the $\G$-orbits in the set of cusps. Fundamental 
results from the theory of automorphic forms assert that: 
\begin{enumerate}
\item[(i)] $\G$ has only finitely many orbits in its set of cusps
and $X^j$ is a compact Hausdorff space.
\item[(ii)] The sheaf of complex valued continuous functions on
$X^j$ that have analytic restrictions to the strata $X$ and $X(j(I))$
gives $X^j$ the structure of a normal analytic variety.
\item[(iii)] The line bundle $\LL^j\to\BB^j$ descends to an analytic 
line bundle on $\G\bs \BB^j$ (a local section is analytic if 
it is continuous and analytic on strata). We denote its sheaf of sections
by $\Lcal$ (the suppression of $j$ in the notation is justified by (v) below).
\item[(iv)] The line bundle $\Lcal$ is ample on $X^\bb$, 
so that the graded algebra 
\[
A^\dot :=\oplus_k H^0(X^\bb, \Lcal^{\otimes k})
\]
is finitely generated with positive degree generators, having $X^{\bb}$ 
as its proj. 
\item[(v)] If $j'$ is an arithmetic system that refines $j$ in the sense
that $j'(I)\subset j(I)$ for all $I$, then the obvious map 
$X^{j'}\to X^j$ is a morphism of analytic spaces and the line bundle $\Lcal$  
on $X^{j'}$ is the pull-back of of its namesake on $X^j$. 
\end{enumerate}
The projective variety $X^{\bb}$ is known as the
\emph{Baily-Borel compactification} and $X^\tor$ as the 
\emph{toric} compactification of $X$. The boundary of $X^{\bb}-X$ 
is finite (its points are called the cusps of $X^\bb$).

Let us see how this works out locally. Given a $k$-isotropic line, then
the orbit space $E_I:=\G_I\bs (\PP(W)-\PP(I^\perp))$ lies as a 
$\CC^\times$-bundle over $X(I)$ and contains $\G_I\bs\BB$ 
as a punctured disc bundle. The obvious morphism $\G_I\bs\BB\to \G\bs\BB=X$ 
restricts to an isomorphism of an open subset of $X$ (in fact, a punctured neighborhood of the cusp defined by $I$) onto a smaller punctured disc bundle
so that the insertion of $X(j(I))$ can be described in terms of 
the $\CC^\times$-bundle  $E_I\to X(I)$. For instance,
the toric compactification amounts to simply filling in the zero section:
we get $E^\tor_I$ as a union of $E_I$ and $X(I)$. The line
bundle $E^\tor_I\to X(I)$ is anti-ample (the real part of $\psi_I$ is a 
Riemann form). This implies that we can contract the zero section 
analytically; the result is a partial Baily-Borel compactification $E^\bb_I$, which adds to $E^I$ a singleton. The intermediate case $E^j_I$ (a union
of $E_I$ and $X(j(I))$) is obtained from $E^\tor_I$ by analytical contraction via the map $X(I)\to X(j(I))$. (This in indeed possible.)

\section{Arithmetic ball arrangements and their automorphic forms} 

Let now be given a collection $\Hcal$ of hyperplanes of $W$ of hyperbolic
signature that is arithmetically defined relative to $\G$.
So $\Hcal$, when viewed as a subset of
$\PP (W^*)$, is a finite union of $\G$-orbits and contained in  
$\PP (W^*)(k)$.
 
We note that for each $H\in\Hcal$, $\BB_H =\PP (H)\cap \BB$ is a symmetric 
subdomain of $\BB$. Its $\G$-stabilizer $\G_H$ is arithmetic in its
$\SU (\psi)$-stabilizer and so the orbit space $X_H:=\G_H\bs \BB_H$ has its own Baily-Borel 
compactification $X_H^{\bb}$.
The inclusion $\BB_H\subset \BB$ induces an analytic morphism 
$X_H\to X$. This morphism is finite and birational
onto its image (which we denote by $D_H$) and extends to a finite morphism $X_H^{\bb}\to X^{\bb}$ whose image is the 
closure of $D_H$ in $X^{\bb}$ (which we denote by $D_H^{\bb}$).
Clearly, $D_H^{\bb}$ only depends on the $\G$-equivalence class of $H$
and so the union $D^{\bb}_\Hcal:=\cup_{H\in\Hcal} D_H^{\bb}$ is a finite one. 
 
\begin{lemma}\label{nonselfinters}
After passing to an arithmetic subgroup of $\G$ of finite index, each
hypersurface $D_H^{\bb}$ is without self-intersection in the sense that
$X_H^{\bb}\to D_H^{\bb}$ is a homeomorphism and $X_H\to D_H$
is an isomorphism.
\end{lemma}
\begin{proof} We do this in two steps.
The set $S_H$ of $\g\in\G$ such that $\g \BB_H\cap \BB_H$ is nonempty
and not equal to $\BB_H$ is a union of double cosets of 
$\G_H$ in $\G$ distinct from $\G_H$. Each such double coset defines 
an ordered pair of distinct irreducible hypersurfaces in $X_H$ with
the same image in $X$ and vice versa. Since $X_H\to D_H$ is
finite and birational, there are only finitely many such pairs of
irreducible hypersurfaces and hence $S_H=\G_H C_H\G_H$ for some finite
set $C_H\subset\G_H$. 
Let $C$ be the union of the sets $C_H$, for which $H$ runs over a system of 
representatives of $\G$ in $\Hcal$. This is a finite subset of 
$\G-\{ 1\}$. If we now pass from $\G$ to a normal arithmetic subgroup of 
$\G$ of finite index which avoids $C$, then the natural maps
$X_H\to X$ are injective and local embeddings. They are also proper, 
and so they are closed embeddings. 

It is still possible that 
$X_H^{\bb}\to X^{\bb}$ fails to be injective for some for some $H$ 
(although this can only happen when $n=2$). We avoid this situation 
essentially  by proceeding as before: given isotropic line 
$I\subset W$ defining a cusp 
then the set $\g\in \G-\G_I$ such that for some $H\in\Hcal$ 
$\g H\not= H$ and $\g H\cap H\supset I$ is a union of double cosets of 
$\G_I$ in $\G$. The number of such double cosets is at most the number of
branches of $D_H^{\bb}$ at this cusp and hence finite. If $C'_I\subset \G-\G_I$ is 
a system of representatives, then let $C'$ be a of union $C'_I$', where 
$I$ runs over a (finite) system of representatives of $\G$ in the set 
$\BB^{\bb}-\BB$. Choose a normal arithmetic subgroup of $\G$ of finite index 
which avoids $C'$. Then this subgroup is as desired.   
\end{proof}

We assume from now on that the hypersurfaces $D_H^{\bb}$ are without 
self-intersection in the sense of Lemma \ref{nonselfinters} and 
we identify $X_H$ with its image $D_H$ in $X$. Then the collection 
$\Hcal_\G$ of the hypersurfaces $D_H^{\bb}$ is an arrangement on 
$X$. We want to determine whether the algebra 
\[
A^\dot_\Hcal :=\oplus_{k=0}^\infty 
H^0(X^\bb, \Lcal(\Hcal_\G)^{(k)}) ) 
\]
is finitely generated and if so, find its proj. 

The hypersurface $D_H^{\bb}$ on $X^{\bb}$ will in general not support a 
Cartier divisor. On the other hand, the closure $D^\tor_H$ of $D_H$ in 
$X^\tor$ is even smooth. So the blow-up of $D_H^{\bb}$ in $X^{\bb}$ will be 
dominated by $X^\tor$. This is then also true for the minimal normal 
blow-up of $X^{\bb}$ for which the strict transform of each $D_H^{\bb}$ is a 
Cartier divisor. It is easy to describe this intermediate blow-up of
$X^{\bb}$. For this, we let for a $k$-isotropic line $I\subset W$,
$j_\Hcal (I)$ be the common intersection of $I^\perp$ and the 
members of $\Hcal$ that contain $I$.   

\begin{lemma}
The assignment $j_\Hcal$ is an arithmetic system relative to $\G$
and the corresponding morphism $X^{j_\Hcal}\to X^{\bb}$
is the smallest normalized blow-up that makes the fractional ideals 
$\Ocal_{X^{\bb}}(D_H^{\bb})$ invertible. The strict transform of 
$D_H^{\bb}$ in $X^{j_\Hcal}$ is an integral Cartier divisor.
\end{lemma}
\begin{proof}
Let $I\subset W$ be a $k$-isotropic line. Then $J:=j_\Hcal (I)$ 
contains $I$, is contained in $I^\perp$ and is defined over $k$ since it is
an intersection of hyperplanes defined over $k$. 
Let $H\in\Hcal$ contain $I$. Then $H$ defines a hypersurface $E_{I,H}^{\bb}$ in 
$E_I^{\bb}$ and an abelian subtorsor $X(I)_H$ of $X(I)$ with quotient 
$X(H\cap I^\perp)$. The latter is a genus one curve which apparently
comes with an origin $o$ (it is an elliptic curve) and 
the complement of this origin of is affine.
Any regular function on the affine curve $X(H\cap I^\perp)-\{ o\}$
with given pole order at $o$ lifts to a regular function on 
$E_I^\tor-E_{I,H}^\tor$ with the same pole order along $E_{I,H}^\tor$ 
(and zero divisor disjoint from the polar set). 
This proves that the effect of the normalized blowing 
up of the fractional ideal $\Ocal_{E_I^{\bb}}(E_{I,H}^{\bb})$ is to replace 
the vertex $X(I^\perp)$ of $E_I^\bb$ by 
$X(H\cap I^\perp)$. The strict transform of $E_{I,H}^{\bb}$ 
in this blow-up is clearly 
a Cartier divisor. If we do this for all the $H\in\Hcal$ that contain $I$, 
then the vertex gets replaced by $X(J)$ as asserted.
\end{proof}

With Convention \ref{convention} in mind we write $X^j$ for 
$X^{j_\Hcal}$.

\begin{proposition}\label{generation} 
For $l$ sufficiently large, $\Lcal (\Hcal_\G)^{(l)}$
is generated by its sections.
\end{proposition} 

We prove this proposition via a number of lemma's.

Let $\Ocal$ be a $\G$-orbit in $W(k)-\{ 0\}$. The projectivization of the hyperplane perpendicular to $w\in\Ocal$ meets $\BB$ precisely when $\psi (w,w)<0$ and in that case
$\psi (\;\; ,w)^{-1}$ defines a meromorphic section of $\Ocal_\BB (-1)$.
Let $l$ be an integer $>2n+1$ and consider the \emph{Poincar\'e-Weierstra\ss}
series
\[
F_\Ocal^{(l)}:=\sum_{w\in\Ocal} \psi (\;\; ,w)^{-l}.
\] 

\begin{lemma}\label{interiorgen} 
Let $K\subset \LL^\times$ be a compact subset. Then the set $\Ocal_0$ of $w\in\Ocal$ for which $\PP(w^\perp)$ meets $K$ is finite and the series
\[
\sum_{w\in\Ocal -\Ocal_0} \psi (\;\; ,w)^{-l}
\]
converges uniformly and absolutely on $K$. In particular,
$F_\Ocal^{(l)}$ represents a meromorphic 
function on $\LL^{\times}$ that is homogeneous of degree $-l$.
\end{lemma}
\begin{proof} 
Since $\Ocal$ is an orbit of the arithmetic group $\G$ in $W(k)-\{ 0\}$,
the $\Ocal_k$-submodule $\Lambda\subset W$ spanned by this orbit is  
discrete in $W$. Now $\Ocal$ 
is also contained in the real hypersurface $\Sigma_a\subset W$ defined by 
$\psi (w,w)=a$ for some $a<0$. The map $m: \Sigma_a\to [0,\infty)$ defined by 
$m(w):=\min_{z\in K} |\psi (z,w)|$ is proper. A standard argument
shows that the cardinality of 
$\Lambda\cap \Sigma_a\cap m^{-1}(N-1,N]$ is proportional to the $(2n+1)$-volume of 
$\Sigma_a\cap m^{-1}(N-1,N]$ and hence is $\le cN^{2n}$ for some $c>0$. 
So $\Ocal_0$ is finite and for $z\in K$ we have the uniform estimate
\[
\sum_{w\in\Ocal -\Ocal_0} |\psi (z,w)|\le \text{constant}+
\sum_{N=2}^\infty cN^{2n}(N-1)^{-l}.
\]
The righthand side converges since $l> 2n+1$.
\end{proof}

For $\Ocal$ and $l$ as above, we extend $F_\Ocal^{(l)}$ to each stratum $\LL^\times (I)$ of $(\LL^\times)^\tor$ by taking 
the subseries defining $F_\Ocal^{(l)}$ whose terms makes sense on that stratum:
\[
F_\Ocal^{(l)}(z):=\sum_{w\in\Ocal\cap I^\perp} \psi (z,w)^{-l} \text{ if } z\in \LL^\times (I).
\]
It is clear from the preceding that this subseries represents a meromorphic
function on $\LL^\times (I)$. The bundle $\LL^\times (I)\to\BB (I)$ is canonically isomorphic to the trivial bundle with fiber $W/I^\perp-\{ 0\}$.
So after choosing a generator of $W/I^\perp$ we can think 
of this subseries as defining a rational function on $\BB (I)$.

\begin{lemma}\label{cuspgen} 
The $\G$-invariant piecewise
rational function on $(\LL^\times)^\tor$ defined by $F_\Ocal^{(l)}$ 
represents a meromorphic section
of $\Lcal$ over the toroidal compactification $X^\tor$ of $X$.
\end{lemma}
\begin{proof}
Given the $k$-isotropic line $I\subset W$, choose a $k$-generator
$e\in I$ and consider the action of the right half plane on $\LL^\times$ 
defined by $T_e$ (see Subsection \ref{stabcusp}):
$T_e^s(z)=z+s\psi (z,e)e$.
In view of the characterization in Subsection \ref{analyticstructure} 
of the analytic structure on $(\LL^\times)^\tor$ it is enough to
show that $\lim_{\text{Re}(s)\to +\infty} F_\Ocal^{(l)}T_e^s |K$
exists as a meromorphic function and this limit is equal to the subseries 
defined above. We claim that with $K$, $\Sigma_a$ and $\Lambda$ as in the 
proof of Lemma \ref{interiorgen} above, the map
\[
(z,w)\in K\times (\Lambda\cap (\Sigma_a\setminus I^\perp))\mapsto
\text{Re}\Big(\frac{\psi (z,w)}{\psi (z,e)\psi (e,w)}\Big)
\]
is bounded from below, say by $\ge -s_o$. Once we establish this we are done,
for then  
\[
s\mapsto |\psi (T_e^sz,w)|=|\psi (z,e)\psi (e,w)|.
\Big\vert\frac{\psi (z,w)}{\psi (z,e)\psi (e,w)}+ \text{Re}(s)\Big\vert
\]
is monotone increasing for $\text{Re}(s)\ge s_o$ for any 
for $z\in K$ and $w\in \Ocal\setminus I^\perp$, implying that the limit is as stated.

In order to prove our claim, we first note that without any loss of generality we may assume that $\psi (z,e)=1$ for all $z\in K$ (just replace 
$z\in K$ by $\psi (z,e)^{-1}z$). Write $e_0$ for $e$ and choose an isotropic $e_1\in W(k)$ such that $\psi (e_0,e_1)=1$. 
Since the orthogonal complement of $\CC e_0+\CC e_1$ is negative definite, we 
denote the restriction of $-\psi$ to this complement by $\la\; ,\;\ra$. 
We write any $w\in W$ as $w_0e_0+w_1e_1+w'$ with $w'\perp (\CC e_0+\CC e_1)$.
Notice that the function $w\in \Lambda\setminus I^\perp\mapsto |w_1|=|\psi (w,e)|$ has a positive minimum $c>0$. Then for $z\in K$ and $w\in \Lambda\cap (\Sigma_a\setminus I^\perp)$ we have the estimate 
\begin{align*}
\text{Re}\Big(\frac{\psi (z,w)}{\psi (e,w)}\Big)
&=\text{Re}\Big(\frac{z_0\bar w_1 +\bar w_0-\la z',w'\ra}{\bar w_1}\Big)\\
&=\text{Re}(z_0+|w_1|^{-2}\bar w_0 w_1-\la z',w'/w_1\ra )\\
&=\text{Re}(z_0)+|w_1|^{-2}(a+\la w',w'\ra )-\la z',w'/w_1\ra )\\
&=\text{Re}(z_0)+|w_1|^{-2}a+ \| w'/w_1-\tfrac{1}{2}z'\|^2
-\tfrac{1}{4}\|z'\|^2\\
&\ge \text{Re}(z_0) -c^{-2}a-\tfrac{1}{4}\|z'\|^2. 
\end{align*} 
The claim follows and with it, the lemma.
\end{proof} 

If $I$ and $\Ocal$ are such that $\Ocal\cap I^\perp$ is a single $\G_I$-orbit,
say of $w_0\in I^\perp$, then if $H_0$ denotes the orthogonal complement
of $w_0$ in $W$, $H_0\cap I^\perp$ is independent of the choice of $w_0$.
The same is true for the abelian divisor $X(I)_{H_0}$ in $X(I)$ and   
the restriction of $F_\Ocal^{(l)}$ to the abelian variety $X(I)$ is in fact
the pull-back of a meromorphic function on the elliptic curve
$X(I)/X(I)_{H_0}=X(H_0\cap I^\perp)$. It has the Weierstra\ss\ form in the sense of Lemma \ref{weierstrass} below.

\begin{lemma}\label{weierstrass}
Let $L\subset\CC$ be a lattice. Then for $k\ge 3$ the series
\[
\wp_k(z):=\sum_{a\in L} (z+a)^{-k}
\]
represents a rational function on the elliptic curve $\CC/L$ and
for any $k_0\ge 3$, the $\wp_k$ with $k\ge k_0$
generate the function field of this curve.
\end{lemma}

We omit the proof of this well-known fact.

\begin{proof}[Proof of \ref{generation}]
It suffices to prove this statement for $\Lcal (\Hcal_\G)^{(k)})$ instead of
$\Lcal (\Hcal_\G)$ for some $k>0$.
We first show that for $k> 2n+1$ the sections of
$\Lcal (\Hcal_\G)^{(k)}$ generate this sheaf over $X$. 
Let $H_0\in\Hcal$ and $z_o\in H_0\cap \LL^\times$.
Choose $w_0\in W(k)$ spanning the orthogonal complement of $H_0$ and
let $\Ocal$ denote the $\G$-orbit of $w_0$. Then according to Lemma \ref{interiorgen} $F_\Ocal^{(k)}$ defines a section of $\Lcal (\Hcal_\G)^{(k)}$.
It is in fact a section of the invertible subsheaf 
$\Lcal (D^j_{H_0})^{(k)}$ on $X^j$ and 
nonzero as such at $z_0$. These subsheafs generate 
$\Lcal (\Hcal_\G)^{(k)}$.

So it remains to verify the generating property over a cusp. 
But this follows from the conjunction of Lemma's \ref{cuspgen}
and \ref{weierstrass}.
\end{proof}

This has the following corollary, which we state as a theorem:

\begin{theorem}\label{mainthm:ball}
The line bundle $\Lcal$ and the collection of strict transforms of the
hypersurfaces $D^\bb_H$ in $X^{j_\Hcal}$
satisfy the conditions of Theorem \ref{singversion} so that the diagram of 
birational morphisms and projective completions of $X^\circ$
\[
X^{\bb}\leftarrow X^{j_\Hcal}\leftarrow \tilde X^\Hcal\to \hat X^\Hcal
\]
is defined. The morphism $\tilde X^\Hcal\to X^{\bb}$ is the 
blowup defined by $\Ocal_{X^\bb}(\Hcal_\G)$ on $X^{\bb}$. The coherent 
pull-back of  $\Lcal (\Hcal_\G)$ to $\tilde X^\Hcal$ is semiample
and defines the contraction $\tilde X^\Hcal\to \hat X^\Hcal$.
\end{theorem}

It is worth noting that the difference $\hat X^\Hcal-X^\circ$ need not be a 
hypersurface: if $\Hcal$ is nonempty, then the complex codimension of this 
boundary is the minimal dimension of a nonempty intersection of members of 
$\Hcal$ in $\BB$. In many interesting examples this is $>1$ and in such cases 
a section of $\Lcal^{\otimes k}$ over $X^\circ$ extends
to a section of $\hat\Lcal (\Hcal_\G)^{\otimes k}$. This gives
the following useful application:

\begin{corollary}\label{maincor:ball}
Let $\Hcal$ be an arithmetic arrangement on a complex ball $\BB$ of dimension 
$\ge 2$ that is arithmetic relative to the arithmetic group $\G$.
Suppose that any intersection 
of members of the arrangement $\Hcal$ that meets $\BB$ 
has dimension $\ge 2$ in $\BB$.
Then the algebra of automorphic forms
\[
\oplus_{k\in\ZZ} H^0(\BB^\circ,\Ocal^\an(\LL)^{\otimes k})^\G
\]
(where $\LL$ is the natural automorphic bundle over $\BB$ and
$\BB^\circ$ is the arrangement complement)
is finitely generated with positive degree generators and its proj
is the modification $\hat X^{\Hcal}$ of $X^\bb$. 
\end{corollary}

\section{Ball arrangements that are defined by automorphic forms}\label{productexp}
The hypotheses of Corollary \ref{maincor:ball} are in a sense opposite 
to those that are needed to ensure that an arithmetically defined ball arrangement
is the zero set of an automorphic form on that ball. This will follow from 
Lemma  \ref{product} and the subsequent remark. Assume
we are in the situation of this section and that we
are given a function $H\in \Hcal\mapsto n_H\in \{ 1,2,\cdots\}$ 
that is constant on the $\G$-orbits. This amounts to giving an effective 
divisor on $X^\bb$ supported by $D^\bb_\Hcal$. Let $I\subset W$ be a $k$-isotropic line
contained in a member of $\Hcal$. Then the collection 
$\Hcal_I$ of $H\in\Hcal$ passing 
through $I$  decomposes into a finite collection of $\G_I$-orbits of 
hyperplanes. Each $H\in\Hcal_I$ defines a hyperplane $(I^\perp/I)_H$ in
$I^\perp/I$ with $\G_I$-equivalent hyperplanes defining the same $\G_I$-orbit.

Suppose now that $\sum_H n_H \BB_H$ is the divisor of an automorphic form.
Then Lemma \ref{product} implies that
\[
\sum_{H\in \G_I\bs\Hcal_I} n_H (I^\perp/I)_H
\]
should represent the Euler class of an ample line bundle. 
This Euler class is in fact known to be the negative of the 
Hermitian form on $I^\perp/I$.
This implies that the intersection of the hyperplanes $(I^\perp/I)_H$
is reduced to the origin. In other words: \emph{the collection of $H\in\Hcal$ 
containing a $k$-isotropic line is either empty or has intersection equal to 
that line.} Although this is formally weaker than the 
above property, this simple requirement turns out be already rather strong.

\begin{question}
Suppose that the effective divisor $\sum_H n_H\BB_H$ is at every cusp
the zero divisor of an automorphic function at that cusp. In other words, suppose that the corresponding Weil divisor on the Baily-Borel compactification is in fact a Cartier divisor. Is it then the divisor of an automorphic form?
More generally, if $H\in\Hcal \mapsto n_H$ is $\ZZ$-valued and $\G$-invariant in such a way that $\sum_H n_H\BB_H$ defines a Cartier divisor on the Baily-Borel 
compactification, is  $\sum_H n_H\BB_H$ then the divisor of a meromorphic automorphic form?
\end{question}

In case such an automorphic form exists, we expect it of course to have a product expansion.

\section{Some applications}
There are a number of concrete examples of moduli spaces to which 
Corollary \ref{maincor:ball} applies. To make the descent from the 
general to the special in Bourbakian style, let us start out from the following 
situation: we are given an integral projective variety $Y$ with ample 
line bundle $\eta$ and a reductive group $G$ acting on the pair 
$(Y,\eta)$. Assume also given a $G$-invariant open-dense subset
$U\subset Y$ consisting of $G$-stable orbits. 
Then $G$ acts properly on $U$ and the orbit space $G\bs U$ exists as 
a quasi-projective orbifold with the restriction $\eta |U$
descending to an orbifold line bundle $G\bs(\eta |U)$ over $G\bs U$. 
To be precise, geometric invariant theory tells us that the algebra of 
invariants $\oplus_{k\in\ZZ} H^0(Y,\eta^{\otimes k})^G$
is finitely generated with positive degree generators. Its proj is
denoted $G\bss Y^\ss$, because a point of this proj
can be interpreted as a minimal $G$-orbit in the  semistable locus 
$Y^\ss\subset Y$. It contains $G\bs U$ as an open dense subvariety. 
The orbifold line bundle $G\bs(\eta |U)$ extends to an orbifold line bundle 
$G\bss \eta$ over $G\bss Y^\ss$ in a way that the space of global 
sections of its $k$th tensor power can be identified with 
$H^0(Y,\eta^{\otimes k})^G$.  

\begin{theorem}\label{gitbb} 
Suppose we are given identification of $G\bs (U, \eta |U)$ 
with a pair coming 
from a ball arrangement $(X^\circ,\Lcal |X^\circ)$, such that
\begin{itemize}
\item[(i)] any nonempty intersection of members of the 
arrangement $\Hcal$ with $\BB$ has dimension $\ge 2$ and  
\item[(ii)] the boundary  $G\bss Y^\ss-G\bs U$ 
is of codimension $\ge 2$ in $G\bss Y^\ss$. 
\end{itemize}
Then this identification determines an isomorphism 
\[
\oplus_{k\in\ZZ} H^0(Y,\eta^{\otimes k})^G\cong
\oplus_{k\in\ZZ} H^0(\BB^\circ,\Ocal^\an (\LL)^{\otimes k})^\G,
\]
in particular, the algebra of automorphic forms is finitely generated
with positive degree generators. Moreover, the isomorphism 
$G\bs U\cong X^\circ$ extends to an isomorphism 
$G\bss Y^\ss\cong \hat X^{\Hcal}$.  
\end{theorem}
\begin{proof}
The codimension assumption implies that any section of  
$G\bs(\eta^{\otimes k}|U)$ extends to $(G\bss\eta^\ss)^{\otimes k}$. 
Since $H^0(\BB^\circ,\Ocal(\LL)^{\otimes k})^\G=
H^0(X,\Lcal^{\otimes k})$, the first assertion follows. This induces an
isomorphism between the underlying proj's and so we obtain an isomorphism
$G\bss Y^\ss\cong \hat X^{\Hcal}$, as stated.
\end{proof}

The identification demanded by the theorem will usually come 
from a period mapping. In such cases the codimension hypotheses 
of are often fulfilled. Let us now be more concrete. 

\subsection{Unitary lattices attached to directed graphs}
The only imaginary quadratic fields we will be concerned with are 
the cyclotomic ones, i.e.,  $\QQ (\zeta_4)$ and $\QQ(\zeta_6)$. 
Their rings of integers
are the Gaussian integers $\ZZ[\zeta_4]$ and the Eisenstein integers 
$\ZZ[\zeta_6]$ respectively. 
It is then convenient to have the following 
notation at our disposal. Let $D$ be a finite graph without loops
and multiple edges and suppose all edges are directed, in other words
$D$ is a finite set $I$ plus a collection of $2$-element subsets of
$I$, each such subset being given as an ordered pair.
Then a Hermitian lattice $\ZZ [\zeta_k]^D$ is defined for $k=4,6$ as follows: 
$\ZZ [\zeta_k]^D$ is the free $\ZZ[\zeta_k]$-module on the set of vertices 
$(r_i)_{i\in I}$ of $D$ and a $\ZZ[\zeta_k]$-valued Hermitian form $\psi$ 
on this module defined by
\begin{equation*}
\psi(r_i,r_j)=
\begin{cases}
|1+\zeta_k|^2=k/2 &\text{ if } j=i,\\
-1-\zeta_k &\text{ if $(i,j)$ is a directed edge},\\
0 &\text{ if $i$ and $j$ are not connected}.
\end{cases}
\end{equation*}
We denote the group of unitary transformations of $\ZZ [\zeta_k]^D$
by $\U (D, \zeta_k)$. This is an arithmetic group in the unitary group of 
the complexification of $\ZZ [\zeta_k]^D$.
The elements in $\ZZ [\zeta_k]^D$ of square norm
$k/2$ are called the \emph{roots} of $\ZZ [\zeta_k]^D$.

If $D$ is a forest, then the isomorphism class 
$\ZZ [\zeta_k]^D$ is independent of the way 
the edges are directed. Usually it is only the isomorphism class that 
matters to us, and so in this case there is no need to specify 
these orientations.

\subsection{The moduli space of quartic curves}
It is well-known fact that the anticanonical 
map of a Del Pezzo surface of degree two realizes that surface as a double
cover of a projective plane ramified along a smooth quartic curve and
that this identifies the coarse moduli space of Del Pezzo surfaces of degree two
and that of smooth quartic curves. The latter is also
the the coarse moduli space of nonhyperelliptic genus three curves. 
The invariant theory of quartic curves is classical (Hilbert-Mumford). 
A period map taking values in a ball quotient has been constructed 
by Kond\=o \cite{kondo1}. Heckman recently observed that the relation between the two is covered by Theorem \ref{gitbb} (unpublished) and that the  
situation is very similar to the case of rational elliptic surfaces 
discussed below, it is only simpler. We briefly review this work. 

Let us first recall the invariant theory of quartic curves. 
Fix a projective plane $P$. Put $H_k:=H^0(P,\Ocal_P)$,
so that $P_k:=\PP(H_k)$ is the space of 
effective degree $k$ divisors on $P$.  We take $Y:=P_4$ and 
$\eta:=\Ocal_Y(1)$ and $G:=\SL (H_1)$. Following \cite{mumford} a point of 
$P_4$ is $G$-stable precisely
if the corresponding divisor is reduced and has cusp 
singularities at worst (i.e., locally formally given by $y^2=x^i$ with 
$i=1,2,3$). It is semistable precisely when it is reduced and has 
tacnodes at worst (like $y^2=x^i$ with $1\le i\le 4$) or is a double 
nonsingular conic. A divisor has a minimal strictly semistable orbit 
precisely when it is sum of two reduced conics, at least one of which is 
nonsingular, which are tangent to each other at two points. If
$(x,y)$ are affine coordinates in $P$, then these orbits 
are represented by the quartics $(xy-\lambda_0)(xy-\lambda_1)$, 
with $[\lambda_0:\lambda_1]\in\PP^1$ and so 
$G\bss Y^\ss-G\bs Y^\st$ is a rational curve. 
The nonsingular quartics are stable and define an open subset 
$Y'\subset Y$.

Kond\=o's period map is defined as follows: given a smooth
quartic curve $Q\subset P$, he considers the $\mu_4$-cover of $P$ 
that is totally ramified along $Q$. This is a 
smooth K3-surface of degree $4$ that clearly double covers the degree two 
Del Pezzo surface attached to $P$. For such surfaces there is
a period map taking values in an arithmetic group quotient of a 
type IV domain of dimension $19$. But the $\mu_4$-symmetry makes it 
actually map in a ball quotient of dimension $6$. To be precise,
consider the Hermitian lattice $\ZZ [\zeta_4]^{E_7}$.
Its signature is $(6,1)$ and so its complexification defines a 
$6$-dimensional ball $\BB$.
We form $X:=\U (E_7,\zeta_4)\bs\BB$. It turns out that $\U (E_7,\zeta_4)$ 
acts transitively on the set of cusps so that $X^\bb$ is topologically
the one-point compactification of $X$.
The roots (i.e., the vectors of square norm $2$) come in two $\U (E_7,\zeta_4)$-equivalence classes:
the one represented by the orthogonal complement of any generator $r_i$, 
denoted $\Hcal_n$, and the rest, $\Hcal_h$, represented by 
$r_i-(1+\zeta_4)r_j$, where $r_i,r_j$ are two generators with $\psi (r_i,r_j)=1+\zeta_4$. 
This defines a hypersurface $D$ in $X$ with two irreducible components
$D_n$ and $D_h$ respectively. The hyperplane sections of $\BB$ of type
$\Hcal_h$ are disjoint and so $D_h$ has no self-intersection.
Kond\=o's theorem states that the period map defines isomorphisms
\[
G\bs Y'\cong X-D,\quad G\bs Y^\st\cong X-D_h.
\]
The last isomorphism is covered  by an isomorphism of orbiline bundles
$G\bss \eta\to \Lcal|X-D_h$ and so Theorem \ref{gitbb} applies with 
$U=Y^\ss$: we find that
the period map extends to an isomorphism 
\[
G\bss Y^\ss\cong \hat X^{\Hcal_h}.
\]

\subsection{A nine dimensional ball quotient}
We begin with a Deligne-Mostow example which involves work of Allcock
\cite{allcock}. Fix a projective line $P$,
put $H_k:=H^0(P,\Ocal_P(k))$, and identify $P_k:=\PP (H_k)$ with the 
$k$-fold symmetric product of $P$, or what amounts to the same, the linear 
system
of effective degree $k$ divisors on $P$.  The group $G:=\SL (H_1)$ acts 
on the pair $(P_k,\Ocal_{P_k}(1))$. Following Hilbert, a degree $k$ divisor
is stable resp.\ semistable if and only if all its multiplicities are 
$< k/2$ resp.\ $\le k/2$. We only have a strictly semistable orbit if $k$
is even $\ge 4$ and in that case there is a unique minimal  
such orbit: the divisors with two distinct points of 
multiplicity
$k/2$. The nonreduced divisors define a \emph{discriminant} hypersurface 
in $P_k$, whose complement we shall denote by $P'_k$.

Now take $k=12$. Given a reduced degree $12$ divisor $D$ on $P$ we can form
the $\mu_6$-cover $C\to P$ that is totally ramified over $D$. The 
Jacobian $J(C)$ of $C$ comes with an action of the covering group $\mu_6$ and
the part where that group acts with a primitive character defines a
quotient abelian variety of dimension $10$. The isomorphism type of this
quotient abelian variety is naturally given by a point in a ball quotient.
To be precise, if $A_{10}$ is the string with $10$ nodes, then  
$\ZZ [\zeta_6]^{A_{10}}$ is nondegenerate of signature $(9,1)$ and so
defines a ball $\BB$ of complex dimension $9$. The group 
$\U(A_{10},\zeta_6)$ is arithmetic and so we may form
\[
X:=\U(A_{10},\zeta_6)\bs \BB.
\] 
It acts transitively on the cusps so that the Baily-Borel 
compactification $X^\bb$ of $X$ is a one-point compactification.
The hyperplanes perpendicular to a root makes up an arrangement 
that is arithmetically defined. Since they are transivily permuted by 
$\U(A_{10},\zeta_6)$, they define an irreducible hypersurface $D$ in $X$.
A theorem of Deligne-Mostow \cite{dm} implies that $X-D$ parametrizes the 
isomorphism types of the quotients of the Jacobians that we encounter in the
situation described above and that the `period map'
\[
G\bs P'_{12}\to X-D
\]
is an isomorphism. It also follows from their work (see also \cite{hl}) 
that this period isomorphism extends to isomorphisms
\[
G\bs P_{12}^\st\cong  X,\quad
G\bss P_{12}^\ss\cong X^\bb,
\]
(with the unique minimal strictly semistable orbit corresponding to the unique
cusp) and that the latter identifies $\Lcal^{\otimes 12}$ with
$\G\bs\Ocal_{P_{12}^\ss}(1)$. (In this situation $\G$ is not neat and
so $\Lcal $ is merely an orbiline bundle on $X^\bb$.)

We can also think of $G\bs P'_{12}$ as the moduli space of
hyperelliptic curves of genus $5$ or as the moduli space of $12$-pointed
smooth rational curves (with the $12$ points not numbered). For either
interpretation this moduli space has a Deligne-Mumford compactification
and these two coincide as varieties. It can be verified that the 
Deligne-Mumford compactification 
is just the minimal normal crossings blowup in the sense of
Subsection \ref{minimalblowup} of $D$ in $X^\bb$. 

\subsection{Miranda's moduli space of rational elliptic surfaces}
In this discussion a rational elliptic surface is always assumed to 
have a section. The automorphisms of such a surface permute the sections 
transivitively, so when we are concerned with isomorphism classes, there is
no need to single out a specific section. Contraction of a section yields
a  Del Pezzo surface of degree one and vice versa, so 
the moduli space of rational elliptic surfaces is the same as the 
moduli space of degree one Del Pezzo surfaces. A rational elliptic surface
can always be represented in Weierstra\ss\ form: 
$y^2=x^3+3f_0(t)x+2f_1(t)$, where $x,y,t$ are the affine coordinates of
a $\PP^2$-bundle over $\PP^1$ and $f_0$ and $f_1$ are rational functions 
of degree $4$ and $6$ repectively. (We normalized the coefficients
as to have the discriminant of the surface take the simple form 
$f_0^3+f_1^2$.) This is Miranda's point of departure
\cite{miranda} for the construction of a compactification by means of 
geometric invariant theory of this moduli space. With the notation
of the previous subsection, let $Y$ be the `orbiprojective' space 
obtained as the orbit space 
of $H_4\oplus H_6-\{ (0,0)\}$ relative to the action of the center 
$\CC^\times\subset\GL (H_1)$. It comes with a natural ample orbifold 
line bundle $\eta$ over $Y$ endowed with an action of $G=\SL (H_1)$. 
It has the property that the pull-back of $\Ocal_{P_{12}}$ under 
the equivariant `discriminant morphism' 
\[
\Delta : Y\to P_{12},\quad (f_0,f_1)\mapsto f_0^3+f_1^2
\]
is equivariantly isomorphic to $\eta^{\otimes 6}$.
The above morphism is finite and birational with image a hypersurface.
The preimage $\Delta^{-1}(P_{12}')$, which clearly parametrizes the 
rational elliptic fibrations over $P$ with reduced discriminant,
maps isomorphically to its image in $P_{12}$ and is contained in the $Y^\st$.
We denote its $G$-orbit space by $\Mcal$. But a stable point of 
$Y$ need not have a stable image in $P_{12}$: Miranda shows that the $G$-orbit space of $Y^\st\cap \Delta^{-1}(P_{12}^\st)$, denoted here by $\Mcal\tilde{}$,  
parametrizes the Miranda-stable rational elliptic surfaces with reduced fibers 
and stable discriminant, in other words, the 
allowed singular fibers have Kodaira type  $I_k$ with $k<6$, $II$, $III$ 
or $IV$, whereas the difference $G\bs Y^\st -\Mcal\tilde{}$ parametrizes 
rational elliptic surfaces with an $I_k$-fiber with $6\le k\le 9$. 
The minimal orbits in $Y^\ss-Y^\st$ parametrize rational elliptic surfaces with two $I^*_0$-fibers or a $I^*_4$-fiber (they make up a rational curve). So the projective variety $\Mcal^M:=G\bss Y^\ss$ still parametrizes distinct 
isomorphism classes of elliptic surfaces. 
We regard this variety as a projective completion of $\Mcal\tilde{}$ and call it the \emph{Miranda compactification}. The boundary $\Mcal^M-\Mcal\tilde{}$ 
is of codimension $5$ in $\Mcal^M$ and
$\Mcal^M-\Mcal$ is a hypersurface in $\Mcal^M$.  
\\

\subsection{A period map for rational elliptic surfaces}
Assigning to a rational elliptic fibration with reduced discriminant its
discriminant defines a morphism 
\[
\Mcal=G\bs Y'\to G\bs P'_{12}\to X-D .
\]
and the pull-back of $\Lcal^{\otimes 12}$ can be identified with 
$G\bs \eta^{\otimes 6}$.
It is proved by Heckman-Looijenga \cite{hl} that this morphism
is a closed embedding with image a deleted hyperball quotient. 
To be precise,
choose a vector in $\ZZ[\zeta_6]^{A_{10}}$ of square norm $6$ (for instance, 
the sum of two perpendicular roots)
and let $\Lambda_o\subset\ZZ[\zeta_6]^{A_{10}}$ be its orthogonal complement. 
It is shown in \emph{op.\ cit.}\ that all such vectors are 
$\U (A_{10},\zeta_6)$-equivalent and for that reason the choice of 
this vector is immaterial for what follows. The $\U (A_{10},\zeta_6)$-stabilizer of $\Lambda_o$ acts on 
$\Lambda_o$ through its full unitary group, which we denote by $\G_o$.
The sublattice $\Lambda_o$ defines a hyperball $\BB_o\subset\BB$,
a ball quotient $X_o:=\BB_{o,\G_o}$ and a natural map $p: X_o\to X$.
The restriction of $\Hcal$ to $\BB_o$ is an arrangement on $\BB_o$ that is arithmetically defined and the corresponding hypersurface $D_o$ on $X_o$ is
just $p^{-1}(D)$. It is proved in \cite{hl} that the period map defines an isomorphism
\[
\Mcal=G\bs Y'\cong X_o -D_o.
\]
It is easily seen to extend to a morphism $\Mcal\tilde{}\to X_o$.
This extension is not surjective. This is related to the fact 
that $D_o$ has
four irreducible components, or equivalently, that the restriction 
of the arrangement $\Hcal$ to $\CC\otimes_\Ocal\Lambda_o$ is not a single $\G_o$-equivalence class, 
but decomposes into four such classes. These four classes are distinguished 
by a numerical invariant (for the definition of which we refer to 
\emph{op.\ cit.}) that takes the values $6,9,15,18$. So 
$\Hcal |\CC\otimes_\Ocal\Lambda_o$ is the disjoint union of the $\G_o$-equivalence classes 
$\Hcal(d)$ with $d$ running over these four
numbers and $\Hcal(d)$ defines an irreducible hypersurface 
$D(d)\subset X_o$. The image of the period map is the complement of 
$D(6)\cup D(9)$. In fact, we have an isomorphism
\[
\Mcal\tilde{}\cong X_o-(D(6)\cup D(9))
\]
with the generic point of 
$D(18)$ resp.\ $D(12)$ describing semistable elliptic fibrations with a 
$I_2$ resp.\ $II$-fiber. This isomorphism is covered by an isomorphism of 
orbifold line bundles: $G\bs \eta^{\otimes 6}$ is isomorphic to the 
pull-back of the $\Lcal^{12}$. It is shown in \cite{hl} that  
any nonempty intersection of hyperplane sections of $\BB_o$ taken from
$\Hcal(6)\cup\Hcal(9)$ has dimension $\ge 5$. 
Since $\Mcal^M-\Mcal\tilde{}$ is of codimension $5$ on $\Mcal^M$,
theorem \ref{gitbb} applies with $U:=Y^\st\cap \Delta^{-1}(P_{12}^\st)$ and we find:

\begin{corollary}
The period map extends to an isomorphism of the Miranda moduli space 
$\Mcal^M$ of
rational elliptic surfaces onto the modification $\hat X_o^{\Hcal(6)\cup\Hcal(9)}$ of $X_o^\bb$ defined by the arrangement $\Hcal(6)\cup\Hcal(9)$. 
\end{corollary}

\subsection{The moduli spaces of cubic surfaces and cubic threefolds}
The previous two examples describe the  moduli spaces of Del Pezzo surfaces of degree one and two as a ball quotient and express its GIT compactification 
as a modified Baily-Borel compactification of the ball quotient.
This can also be done for Del Pezzo surfaces of degree three, that is
for smooth cubic surfaces. Allcock, Carlson and Toledo gave in \cite{act} 
a $4$-ball quotient description of this moduli space, by assiging to
a smooth cubic surface $S\subset \PP^3$ first the $\mu_3$-cover 
$K\to\PP^3$ totally ramified over $S$ and then the intermediate Jacobian
of $K$ with its $\mu_3$-action. But in this case they 
find that the GIT compactification \emph{coincides} with the Baily-Borel compactification (one might say that the relevant arrangement is empty).
However, they recently extended this to
the moduli space of smooth cubic threefolds: the GIT compactification
has been recently given by Allcock \cite{allcock2}, whereas
Allcock, Carlson and Toledo found a $10$-ball quotient description of the 
moduli space (yet unpublished). It comes naturally with an arrangement and
it seems that Theorem \ref{gitbb} applies here, too.

\end{document}